\documentclass[a4paper]{amsart} 
\usepackage{amsmath}
\usepackage{epsfig}
\usepackage{amssymb,showkeys} 
\usepackage{xspace}
\usepackage{cite}
\usepackage{float}
\usepackage{algorithm}
\usepackage{algpseudocode}

\newtheorem{theorem}{Theorem}[section]

\newtheorem{definition}[theorem]{Definition}
\newtheorem{remark}[theorem]{Remark}

\newcommand{\I}{\mathrm{i}}

\newcommand{\C}{\mathbb{C}}
\newcommand{\Z}{\mathbb{Z}}

\renewcommand{\P}{\mathbb{P}}

\newcommand{\id}{\mathbb{I}}

\newcommand{\Zg}{\mathbb{Z}^g}
\newcommand{\Cg}{\mathbb{C}^g}
\newcommand{\Rm}{\mathbb{B}}   
\newcommand{\siegel}{\mathbb{H}^g}

\newcommand{\Rs}{\mathcal{R}}
\newcommand{\Av}{\mathcal{A}}		
\newcommand{\abel}{\alpha}

\newcommand{\Jg}{\mathcal{J}_g}
\newcommand{\SpZ}{\text{Sp}(2g,\Z)}
\newcommand{\Kum}{{\rm Kum}}

\newcommand{\vz}{\mathrm{z}}

\newcommand{\vn}{\mathrm{N}}

\newcommand{\B}{\mathbb{B}}

\begin{document}

\title{Computational approach to the Schottky problem}

\author{Eddy Brandon de Leon}
    \email{eddybrandon11@hotmail.com}
\address{Institut de Math\'ematiques de Bourgogne,
		Universit\'e de Bourgogne, 9 avenue Alain Savary, 21078 Dijon
		Cedex, France}
		
\author{J\"org Frauendiener }
\email{joerg.frauendiener@otago.ac.nz}
\address{Department of Mathematics and Statistics, 
University of Otago,      
P.O. Box 56, Dunedin 9010, New Zealand}
\author{Christian Klein}
    \email{Christian.Klein@u-bourgogne.fr}
\address{Institut de Math\'ematiques de Bourgogne, Institut 
Universitaire de France, 
		Universit\'e de Bourgogne, 9 avenue Alain Savary, 21078 Dijon
		Cedex, France}

              \date{\today} 
              \thanks{This work was partially supported by the EIPHI Graduate School (contract ANR-17-EURE-0002) and by the European Union Horizon 2020 research and innovation program under the Marie Sklodowska-Curie RISE 2017 grant agreement no. 778010 IPaDEGAN. JF thanks for the support by the CNRS for an extended stay in Dijon and the hospitality of the IMB. We thank S. Grushevsky, D. Korotkin and B. Sturmfels for helpful discussions and hints.}

\begin{abstract}
  We present a computational approach to the classical Schottky problem based on Fay's trisecant identity for genus $g\geq 4$. For a given Riemann matrix $\mathbb{B}\in\mathbb{H}^{g}$, the Fay identity establishes linear dependence of secants in the Kummer variety if and only if the Riemann matrix corresponds to a Jacobian variety as shown by Krichever. The theta functions in terms of which these secants are expressed depend on the Abel maps of four arbitrary points on
a Riemann surface. However, there is no concept of an Abel map for 
general $\mathbb{B} \in \mathbb{H}^{g}$. To establish linear 
dependence of the secants, four components of the vectors entering 
the theta functions can be chosen freely. The remaining components 
are determined by a Newton iteration to minimize the residual of the 
Fay identity.  Krichever's theorem assures that if this residual 
vanishes within the finite numerical precision for a generic choice 
of input data, then the Riemann matrix is with this numerical 
precision the period matrix of a Riemann surface. The algorithm is 
compared in genus 4 for some examples to the Schottky-Igusa modular 
form, known to give the Jacobi locus in this case. It is shown that 
the same residuals are achieved by the Schottky-Igusa form and the 
approach based on the Fay identity in this case. In genera 5, 6 and 7, we discuss known examples of Riemann matrices and perturbations thereof for which the Fay identity is not satisfied. 
\end{abstract}



\maketitle

\section{Introduction}

The classical Schottky problem is concerned with identifying Riemann 
matrices as period matrices of some compact Riemann surface of genus 
$g\geq4$ among all symmetric matrices with positive definite 
imaginary part forming the Siegel halfspace $\mathbb{H}^{g}$. In this 
paper, we present a computational approach to this problem based on 
Fay's trisecant identity \cite{fay}, which is an identity between 
theta functions depending on the Abel maps of four arbitrary distinct 
points on a Riemann surface. As shown by Krichever \cite{Kri}, this 
identity only holds on Jacobians and thus gives a sharp criterion to 
identify the \textit{Jacobi locus} in $\mathbb{H}^{g}$, namely, the 
subset whose elements are period matrices of some suitable Riemann 
surface of genus $g$. The idea is to decide for a given Riemann 
matrix $\mathbb{B}\in\mathbb{H}^{g}$, and having no information about 
the Abel map, whether it is possible to choose the arguments of the 
theta functions, which are $g$-dimensional vectors, such that the 
trisecant identity is satisfied. This is equivalent to the following 
optimization problem: can the  vectors, corresponding to Abel maps if 
$\mathbb{B}$ is in the Jacobi locus, be selected in a way such that 
the residual of the Fay identity vanishes with numerical precision 
$\delta$\footnote{In principle, the algorithm works for any chosen 
numerical accuracy. We work here with double precision, i.e., a 
maximal accuracy of the order $\delta\sim10^{-16}$; in practice, due to rounding errors, a maximal precision of $\delta\sim10^{-12}$ can be reached in this setting.}? This would indicate that $\mathbb{B}$ is in the Jacobi locus with precision $\delta$ and that the vectors for which Fay's identity is satisfied define Abel maps.

The Schottky problem essentially goes back to Riemann; see 
\cite{RF,Gru,FGSM} for additional references and a more detailed 
history of the problem. The Siegel halfspace $\mathbb{H}^{g}$ of 
Riemann matrices has (complex) dimension $g(g+1)/2$. But the Jacobi 
locus of Riemann matrices, being period matrices of Riemann surfaces 
of genus $g$, has the (complex) dimension $3g-3$. Thus the Schottky problem to 
identify such matrices is non-trivial for $g\geq4$. Schottky solved the problem in genus 4 via what is now known as the 
Schottky-Igusa modular form \cite{Sch} (a rigorous proof was given 
by Igusa in \cite{Igu}), 
a modular invariant built from \emph{theta constants}, theta functions with 
argument equal to zero. The Schottky-Igusa form vanishes exactly on the 
Jacobi locus. A generalization of this form to a higher genus is not 
known, but see for instance \cite{AC} for attempts in this direction. 
Farkas, Grushevsky and Salvati 
Manni gave a solution to the weak Schottky problem in terms of 
modular invariants for arbitrary genus in \cite{FGSM}. This means that 
their modular invariant vanishes on the Jacobi locus, but as was 
shown for genus 5 by Donagi \cite{Don}, vanishing of the invariant 
for some matrix $\mathbb{B}$
does not guarantee this matrix is related to a Jacobian. 

Algebro-geometric approaches to integrable systems, see for instance 
\cite{algebro,Dintro}, not only allowed to construct quasi-periodic 
solutions to integrable partial differential equations (PDEs), but 
also proved useful in the context of the Schottly problem. 
Krichever's solution of the Kadomtsev-Petviashvili (KP) equation in terms 
of multi-dimensional theta functions on arbitrary compact Riemann 
surfaces raised the question of the computability of the solutions. An 
effective parametrization of the solutions will not work for $g\geq4$ 
exactly because of the Schottky problem. This led Novikov to 
conjecture that Krichever's formula gives a solution to the KP 
equation if and only if the matrix $\mathbb{B}$ is in the Jacobi locus. Shiota 
proved this conjecture in \cite{Shi}. Since it was stressed in 
\cite{tata} that many integrable PDEs can be solved in terms of theta 
functions as a consequence of Fay's identity (\ref{fayid}), Welters 
\cite{Wel} conjectured that this identity holds as the Novikov 
conjecture just on Jacobians. This was proven by Krichever \cite{Kri} 
and is the basis of the algorithm to be presented in this paper. 

To state Fay's identity, we consider four arbitary points on a 
Riemann surface $\mathcal{R}$ of genus $g$:   
$P_{1},P_{2},P_{3},P_{4}\in\mathcal{R}$. The \emph{Abel map}
\begin{align*}
\abel : \Rs & \to Jac(\Rs),\\
 P & \mapsto \left(\int_{P_0}^P \omega_1, ... , \int_{P_0}^P \omega_g \right) \mod \Lambda,
\end{align*}
is a bijective map from the surface into the
Jacobian $Jac(\mathcal{R}) = \mathbb{C}^{g}/\Lambda$, where $\Lambda$ is
the full-rank lattice formed by the periods of the holomorphic 
one-forms on $\Rs$. Moreover, we consider a normalized basis 
$\{\omega_1,...,\omega_g\}$ of differentials with respect to a 
symplectic basis of the first homology group 
$\{a_1,b_1,...,a_g,b_g\}$, i.e., $\int_{a_j}\omega_k = \delta_{jk}$. 
With this choice, the matrix of $b$-periods, which has components $\B_{jk} = \int_{b_j} \omega_k$, is symmetric and with positive definite imaginary part. Thus, the lattice $\Lambda$ takes the form 
\[
\Lambda = \Zg + \B \Zg =
\left\{\lambda\in\C^g \quad | \quad \lambda = \mathrm{m} +\mathbb{B} \,\mathrm{n}; \quad \mathrm{m}, \mathrm{n}\in \mathbb{Z}^{g}\right\}.
\]
The Abel map is extended to any divisor $D=\sum_{k=1}^{N} n_k P_{k}$ on $\Rs$ 
(a divisor is  a formal symbol, $n_{k}\in\mathbb{Z}$, 
$k=1,\ldots,N$) via $\abel(\sum_{k=1}^{N} n_k P_{k})=\sum_{k=1}^{N} 
n_k \abel(P_{k})$.  Following \cite{FaKr} we introduce the cross
ratio function ($\Theta^{*}$ is a theta function with a nonsingular 
odd
characteristic, see (\ref{theta}))
\begin{equation}
    \lambda_{1234}=\frac{\Theta^{*}(\abel(P_{1})-\abel(P_{2}))\Theta^{*}(\abel(P_{3})-\abel(P_{4}))}{
    \Theta^{*}(\abel(P_{1})-\abel(P_{4}))\Theta^{*}(\abel(P_{3})-\abel(P_{2}))}
    \label{eq:cross3}\;,
\end{equation}
which is a function on $\mathcal{R}^4$ that vanishes for $P_{1}=P_{2}$
and $P_{3}=P_{4}$ and has poles for $P_{1}=P_{4}$ and $P_{2}=P_{3}$.
Then the following
    identity holds:
    \begin{equation}
	\begin{split}
    0=&\lambda_{3124}\,\Theta(\mathrm{u}+\smallint_{P_{2}}^{P_{3}})\,
                \Theta(\mathrm{u}+\smallint_{P_{1}}^{P_{4}}) +\lambda_{3214}\,
        \Theta(\mathrm{u}+ \smallint_{P_{1}}^{P_{3}})\,\Theta(\mathrm{u}+
                \smallint_{P_{2}}^{P_{4}})\\
        &-\Theta(\mathrm{u})\;\Theta(\mathrm{u}+\smallint_{P_{2}}^{P_{3}}+
        \smallint_{P_{1}}^{P_{4}})\;,
    \label{fayid}
    \end{split}
\end{equation} 
for all $\mathrm{u}\in \mathbb{C}^{g}$, where we use the notation $\int_{P}^{Q}=\abel(Q-P)$ as well as $\int_{P_1+P_2}^{Q_1+Q_2}=\int_{P_1}^{Q_1}+\int_{P_2}^{Q_2}=\abel(Q_1+Q_2-P_1-P_2)$. 

With the binary addition theorem (\ref{binary}) for theta functions,
we get for (\ref{fayid})
\begin{equation}
	\lambda_{3124}\Theta\left[\epsilon/2\atop 0 
	\right]\left(\int_{P_2+P_4}^{P_1+P_3},2\mathbb{B}\right)
	+\lambda_{3214}\Theta\left[\epsilon/2\atop 0 \right]\left(\int_{P_1+P_4}^{P_2+P_3},2\mathbb{B}\right)
	-\Theta\left[\epsilon/2\atop 0 
	\right]\left(\int_{P_1+P_2}^{P_3+P_4},2\mathbb{B}\right)=0
	\label{Fay2},
\end{equation}
for all $\epsilon \in \Z^g/2\Z^g $, i.e., the $g$ dimensional vectors $\epsilon$ with components 
$\epsilon_{i}=0,1$, $i=1,\ldots,g$. Thus the system (\ref{Fay2}) 
consists of 
$2^{g}$ equations. 
We put 
\begin{equation}
	X := \frac{1}{2} \int_{P_1+P_2}^{P_3+P_4},\quad 	
	Y := \frac{1}{2} \int_{P_2+P_3}^{P_1+P_4},\quad 	
	Z  : = \frac{1}{2} \int_{P_2+P_4}^{P_1+P_3},
	\label{XYZ}
\end{equation}
which means that the $2^{g}$ equations (\ref{Fay2}) take the form 
\begin{equation}
	F_{\epsilon}(X,Y,Z):=c_{1}\Theta\big[\begin{smallmatrix} \epsilon/2 \\ 0 \end{smallmatrix} \big] (2Z,2\mathbb{B})+c_{2}\Theta\big[\begin{smallmatrix} \epsilon/2 \\ 0 \end{smallmatrix} \big] (2Y,2\mathbb{B})-\Theta \big[\begin{smallmatrix} \epsilon/2 \\ 0 \end{smallmatrix} \big] (2X,2\mathbb{B})
	\label{Fay3},
\end{equation}
where 
\begin{equation}
	c_{1} = \frac{\Theta^{*}(Y-X)\Theta^{*}(Y+X)}{
	\Theta^{*}(Y+Z)\Theta^{*}(Y-Z)},\quad 
	c_{2} = \frac{\Theta^{*}(X+Z)\Theta^{*}(X-Z)}{
	\Theta^{*}(Y+Z)\Theta^{*}(Y-Z)}
	\label{c}.
\end{equation}

In this paper we use equation \eqref{Fay2} in the following way: we 
choose $X$, $Y$, $Z$ to be in the fundamental domain given by
\begin{equation}
U_{\B} := \left\{ \mathrm{z}\in \C^g \quad | \quad
\mathrm{z} = \mathrm{p} + \mathbb{B}\,\mathrm{q},\quad p_{i},q_{i} \in [-1/2,1/2], \quad i=1,\ldots,g \right\}.	
	\label{fund}
\end{equation}
In other words,  an arbitrary point of the Jacobian 
can be given in terms of the $(p_{i}, q_{i})$, which are called 
\emph{characteristics}. Thus, the function $F$ with components \eqref{Fay3} 
has the form $F:U^3_{\B}\subset \C^{3g} \to \C^{2^g}$, where $U^3_{\B}$ is a compact subset, and it is known from Fay's identity, see \cite{tata}, that if $\C^g/\Lambda$ is a Jacobian the set of 
zeros of $F$ is 4-dimensional, and it is parametrized according to \eqref{XYZ}. However, we are interested in general $\C^g/\Lambda$ where the notion of an Abel map is unknown; therefore, the zeros of $F$ cannot be determined through \eqref{XYZ}. Thus, if we fix for $X$, $Y$, $Z$ the 
first components and for instance the second component of $X$,  the function $F$ 
can be seen as a function of $\mathbf{x}$, the 
vector built from the remaining components of $X$, $Y$, $Z$ of 
dimension $3g-4$. According to 
Krichever's theorem, this function 
vanishes for generic values of $X_{1,2}$, $Y_{1}$ and $Z_{1}$ for a 
certain choice of the vector 
$\mathbf{x}$ if  and only if
$\mathbb{B}$ is in the Jacobi locus. In the context of the 
Schottky problem, the task is to find the zero of $F(X,Y,Z)$ in the 
space of dimension $3g-4$ by varying $\mathbf{x}$. If the residual is 
 smaller than the 
numerical accuracy $\delta$, then the matrix is within numerical 
precision in the Jacobi locus. Since $F$ is a locally holomorphic function 
of $X$, $Y$, $Z$, it is possible to find the zero with a standard 
Newton iteration as will be shown for examples in genus $g\leq 7$.

The paper is organized as follows: In section~2, we summarize basic 
definitions regarding theta functions, the Siegel fundamental domain and the 
Schottky-Igusa form. In section~3 we present the numerical 
algorithms. In section~4, we consider known examples in genus 4 and 
compare our algorithm to the Schottky-Igusa form. We show that both 
the Schottky-Igusa form and the approach based on the Fay identity 
yield the same residuals for the studied examples. This means they 
identify the Jacobi locus with the same numerical accuracy.  In
section~5, we consider examples in genera 5, 6 and 7. 
We add some concluding remarks in section~6. A pseudo code for the 
algorithm is given in the appendix.

\section{Preliminaries}
In this section, we summarize some theoretical facts needed in the 
following to numerically address the Schottky problem. 

\subsection{Jacobians and principally polarized Abelian varieties}
Let $\siegel$ be the space of $g\times g$ complex symmetric matrices 
with positive definite imaginary part called the Siegel halfspace. The Riemann theta function is given by the series
\begin{equation} \label{theta0}
    \Theta(\vz,\Rm) = \sum_{\vn\in\Z^g} \exp( \pi \I \langle \vn, \Rm \vn \rangle + 2\pi \I \langle \vn,\vz \rangle )
\end{equation}
with $\vz \in \mathbb{C}^{g}$ and where $\left \langle\cdot,\cdot\right\rangle$
denotes the Euclidean scalar product $\left\langle \vn,\vz \right\rangle=\sum_{i=1}^g \vn_i \vz_i$. Since $\Rm \in \siegel $, the series \eqref{theta0} converges uniformly for all $\vz\in\C^g$; thus it is an entire function. \\

The Siegel space parametrizes the set of principally polarized 
Abelian varieties (PPAV) via the assignment $\Rm \mapsto 
(\Av_{\Rm},\Theta_\Rm)$, where the complex torus has the form 
$\Av_{\Rm}=\C^g/(\Z^g+\Rm\Zg)$ and $\Theta_\Rm$ is the divisor of the 
theta function \eqref{theta0} (the set of zeros of the theta 
function), which is a well-defined subvariety of $\Av_{\Rm}$ because of its quasi-periodicity properties \eqref{eq:periodicity}. In the following, we refer to the Riemann matrix $\Rm$ and the PPAV $(\Av_{\Rm},\Theta_\Rm)$ interchangeably. \\

In the special case where $\B$ is the period matrix of some Riemann surface $\Rs$ (which is only unique up to a modular transformation), its associated PPAV is known as the Jacobian $Jac(\Rs)$ of $\Rs$. 

\subsection{Symplectic transformations} 
Two Riemann matrices $\Rm,\tilde{\Rm}\in\siegel$ can define 
isomorphic PPAVs if they define equivalent complex tori and if they are principally polarized. $\Av_{\tilde{\Rm}}$ and $\Av_{\Rm}$ are equivalent as complex tori if there exists an invertible homomorphism $\Cg/\tilde{\Lambda} \to \Cg/\Lambda$. This is equivalent to the existence of a linear transformation $\mathcal{M}: \Cg\to\Cg$ with $\mathcal{M}(\tilde{\Lambda})=\Lambda$ and a matrix
\begin{equation} \label{siegel1}
R=\begin{pmatrix}
A & B\\
C & D
\end{pmatrix} \in M_{2g}(\Z)
\end{equation}
such that
\begin{align} \label{action_Sp}
\mathcal{M} (\tilde{\Rm},\id_g) = (\Rm,\id_g) R^\top,
\end{align}
holds. Additionally, in order for the principal polarization to be preserved, the condition
\begin{align} \label{siegel2}
\begin{pmatrix}
A & B\\
C & D
\end{pmatrix}^\top \begin{pmatrix}
0_g & \id_g\\
-\id_g & 0_g
\end{pmatrix} \begin{pmatrix}
A & B\\
C & D
\end{pmatrix} = \begin{pmatrix}
0_g & \id_g\\
-\id_g & 0_g
\end{pmatrix}
\end{align}
must be satisfied \cite{BL}, which means that $R\in \SpZ$. Thus, two matrices in $\siegel$ conjugated under the action of the symplectic (also called modular) group $\SpZ$ define equivalent PPAVs, where the action of $\SpZ$ on $\siegel$ is given by \eqref{action_Sp}, which is equivalent to
\begin{align} \label{siegel3}
\tilde{\Rm} = R \cdot \Rm = (A\Rm + B)(C \Rm + D)^{-1},
\end{align}
which we call modular transformation. Thus, given $\Rm\in\siegel$ we can choose the most convenient symplectically equivalent $\tilde{\Rm}$ in order to perform the computations as efficiently as possible. For this we will consider Siegel's fundamental domain.\\

Siegel~\cite{Sie} gave the following fundamental domain for the modular group:
\begin{definition}\label{def}
 Siegel's fundamental domain is the subset of $\mathbb{H}^{g}$   such that 
 $\mathbb{B}=X+\mathrm{i}Y\in\mathbb{H}^{g}$ satisfies:
 \begin{enumerate}
     \item  $|X_{nm}|\leq 1/2$, $n,m=1,\ldots,g$,
 
     \item  $Y$ is in the fundamental region of Minkowski reductions 
     \cite{Min1,Min2},
 
     \item  $|\det(C\mathbb{B}+D)|\geq 1$ for all $C$, $D$  (\ref{siegel2}).
 \end{enumerate}
\end{definition}

\subsection{Theta functions}
The multi-dimensional theta functions with characteristics are defined as
the series, 
\begin{equation}\label{theta}
    \Theta \big[\begin{smallmatrix} \mathrm{p} \\ \mathrm{q} \end{smallmatrix} \big](\mathrm{z},\mathbb{B}) = \sum\limits_{\mathrm{N}\in\mathbb{Z}^g}\exp\left\{
    \pi\I\left\langle \mathrm{N}+\mathrm{p},\mathbb{B}\left(\mathrm{N}+\mathrm{p}\right)
    \right\rangle+2\pi \I
    \left\langle \mathrm{N}+\mathrm{p}, \mathrm{z}+\mathrm{q}
    \right\rangle\right\}
    \;,
\end{equation}
with the \emph{characteristics} $\mathrm{p}$,
$\mathrm{q}\in{ \mathbb{R}}^g$. Analogously to \eqref{theta0}, the theta function with characteristics is an entire function on $\C^g$. If $\mathbb{B}$ is the period matrix of some Riemann surface $\Rs$, a characteristic is called \emph{singular} if the corresponding theta function vanishes identically on $\Rs$.

Of special interest are half-integer characteristics with $2\mathrm{p},2\mathrm{q}\in \mathbb{Z}^{g}$.  
Such a characteristic is called \emph{even} if $4\langle 
\mathrm{p},\mathrm{q}\rangle=0\mbox{ mod } 2$ and \emph{odd} 
otherwise. It can be easily shown that 
theta functions with odd (even) characteristics are odd
(even) functions of the argument $\mathrm{z}$.  The theta function with
characteristic is related to the Riemann theta function $\Theta$, the
theta function with zero characteristic $\Theta:= \Theta \big[\begin{smallmatrix} 0 \\ 0 \end{smallmatrix} \big]$,
via
\begin{equation}
    \Theta \big[\begin{smallmatrix} \mathrm{p} \\ \mathrm{q} \end{smallmatrix} \big](\mathrm{z},\mathbb{B}) = \Theta(\mathrm{z}
    +\mathbb{B}\mathrm{p} + \mathrm{q})\exp\left\{\pi\I
    \left\langle\mathrm{p}, \mathbb{B}\mathrm{p}\right\rangle+
    2\pi \I\left\langle\mathrm{p},\mathrm{z} + \mathrm{q}\right\rangle
    \right\}\;.
    \label{thchar}
\end{equation}
Note that from now on we will suppress 
the second argument of the theta functions  if it is 
just $\mathbb{B}$ for ease of 
representation. 
From its definition, a theta function has the periodicity properties 
\begin{equation}
    \Theta \big[\begin{smallmatrix} \mathrm{p} \\ \mathrm{q} \end{smallmatrix} \big](\mathrm{z} +\mathrm{e}_{j}) = 
    \mathrm{e}^{2\pi \mathrm{i}p_{j}}
    \Theta \big[\begin{smallmatrix} \mathrm{p} \\ \mathrm{q} \end{smallmatrix} \big] (\mathrm{z})\;,
    \quad 
    \Theta \big[\begin{smallmatrix} \mathrm{p} \\ \mathrm{q} \end{smallmatrix} \big](\mathrm{z} +\mathbb{B}
    \mathrm{e}_{j})=
    \mathrm{e}^{-2\pi \I (z_{j}+q_{j}) - \I\pi B_{jj}}
    \Theta \big[\begin{smallmatrix} \mathrm{p} \\ \mathrm{q} \end{smallmatrix} \big](\mathrm{z})\;
    \label{eq:periodicity},
\end{equation}
where $\mathrm{e}_{j}$ is a vector in $\mathbb{Z}^{g}$ consisting of zeros
except for a $1$ in the $j$-th position. 
Theta functions satisfy the \emph{binary addition theorem}, see for 
instance \cite{Tai}
\begin{equation}
	\begin{split}
		&\Theta\big[\begin{smallmatrix} \alpha \\ \gamma \end{smallmatrix} \big](z_{1}+z_{2},\mathbb{B})
	\Theta\big[\begin{smallmatrix} \beta \\ \delta \end{smallmatrix} \big](z_{1}-z_{2},\mathbb{B})=\\
	&\sum_{\epsilon\in\mathbb{Z}^{g}/2\mathbb{Z}^{g}}^{}\Theta\big[\begin{smallmatrix}(\alpha+\beta)/2+\epsilon/2 \\ \gamma+\delta \end{smallmatrix} \big] (2z_{1},2\mathbb{B})
	\Theta \big[\begin{smallmatrix} (\alpha-\beta)/2+\epsilon/2 \\
	\gamma-\delta \end{smallmatrix} \big](2z_{2},2\mathbb{B}).
	\end{split}
	\label{binary}
\end{equation}

The action of the modular group on theta functions is known, see for 
instance \cite{algebro,fay,tata}. One has
\begin{equation}
    \Theta \big[\begin{smallmatrix} \tilde{\mathrm{p}} \\ \tilde{\mathrm{q}} \end{smallmatrix} \big]
(\mathcal{M}^{-1}\mathrm{z},\tilde{\mathbb{B}}) = 
    k\sqrt{\det(\mathcal{M})}\exp\left(\frac{1}{2}\sum_{i\leq 
    j}^{}z_{i}z_{j}\frac{\partial}{\partial 
    \mathbb{B}_{ij}}\ln\det\mathcal{M}\right)\Theta \big[\begin{smallmatrix} \mathrm{p} \\ \mathrm{q} \end{smallmatrix} \big] (\mathrm{z}),
    \label{thetamod}
\end{equation}
where $\tilde{\mathbb{B}}$ is given by (\ref{siegel3}) and $k$ is 
a constant with respect to $\mathrm{z}$, and where
\begin{equation}
    \mathcal{M}=C\mathbb{B}+D,\quad 
    \begin{pmatrix}
        \tilde{\mathrm{p}} \\
        \tilde{\mathrm{q}}
    \end{pmatrix}=
    \begin{pmatrix}
        D& -C \\
        -B & A
    \end{pmatrix}    \begin{pmatrix}
        \mathrm{p} \\
        \mathrm{q}
    \end{pmatrix}
+\frac{1}{2}  \begin{pmatrix}  \mbox{diag}(CD^\top) \\
       \mbox{diag} (AB^\top)
    \end{pmatrix}
    \label{thetamod2},
\end{equation}
with diag denoting the diagonal of the corresponding matrices.

%

\subsection{Schottky-Igusa modular form}
The Schottky-Igusa modular form \cite{Sch,Igu} is a polynomial of 
degree 16 in the theta constants $\Theta \big[\begin{smallmatrix} \mathrm{p} \\ \mathrm{q} \end{smallmatrix} \big] (0)$. 
We follow here the presentation in \cite{CKS}: choose the 
characteristics
\begin{equation}
    \mathrm{p}^{(1)}=\mathrm{q}^{(1)}=\frac{1}{2}
    \begin{pmatrix}
        1 \\
        0 \\
        1 \\
        0
    \end{pmatrix},\quad 
     \mathrm{p}^{(2)}=\frac{1}{2}
    \begin{pmatrix}
        0 \\
        0 \\
        0 \\
        1
    \end{pmatrix},
    \quad 
     \mathrm{q}^{(2)}=\frac{1}{2}
    \begin{pmatrix}
        1 \\
        0 \\
        0 \\
        0
    \end{pmatrix},
    \nonumber
\end{equation}
and 
\begin{equation}
    \mathrm{p}^{(3)}=\frac{1}{2}
    \begin{pmatrix}
        0 \\
        0 \\
        1 \\
        1
    \end{pmatrix},\quad 
     \mathrm{q}^{(3)}=\frac{1}{2}
    \begin{pmatrix}
        1 \\
        0 \\
        1 \\
        1
    \end{pmatrix}.
    \nonumber
\end{equation}
Consider  a rank 3 subgroup $N$ of $\mathbb{Z}^8/(2\mathbb{Z}^{8})$ 
generated by 
\begin{equation}
    \mathrm{n}_{1}=\frac{1}{2}
    \begin{pmatrix}
        0 &1\\
        0 &1\\
        0 &1\\
        1 & 0
    \end{pmatrix},\quad  
        \mathrm{n}_{2}=\frac{1}{2}
    \begin{pmatrix}
        0 &0\\
        0 &0\\
        1 &0\\
        1 & 1
    \end{pmatrix},\quad  
        \mathrm{n}_{3}=\frac{1}{2}
    \begin{pmatrix}
        0 &1\\
        0 &0\\
        1 &1\\
        0 & 1
    \end{pmatrix}.
    \nonumber
\end{equation}
We consider the product of theta constants
\begin{equation}
    \pi_{i} = 
    \prod_{(\mathrm{p}\mathrm{q})\in(\mathrm{p}^{(i)}\mathrm{q}^{(i)})+N}^{}\Theta \big[\begin{smallmatrix} \mathrm{p} \\ \mathrm{q} \end{smallmatrix} \big] (0)
    \label{pi}
\end{equation}
and define the Schottky-Igusa modular form as
\begin{equation}
    \Sigma(\Rm) = 
    \pi_{1}^{2}+\pi_{2}^{2}+\pi^{2}_{3}-\pi_{1}\pi_{2}-\pi_{1}\pi_{3}-\pi_{2}\pi_{3}
    \label{SI}.
\end{equation}
It was shown in \cite{Sch,Igu} that this modular form vanishes 
exactly on the Jacobi locus. Our numerical approach is equivalent to \eqref{SI} in the sense that the only input we require is the matrix $\Rm\in\siegel$.

\subsection{Kummer variety}
For the geometrical interpretation of Fay's identity, we need to 
introduce the Kummer variety, which is the embedding of the complex 
torus into projective space via level-two theta functions (theta 
functions of double period). 
\begin{definition}
Let $(\Av_\Rm,\Theta_\Rm)$ be an indecomposable PPAV. Thus, its Kummer variety is the image of the map
\begin{align*}
\Kum: \Av_\Rm/\sigma & \to \P^{2^g-1},\\
Z & \mapsto  \left[\Theta \big[\begin{smallmatrix} \epsilon/2 \\ 0 \end{smallmatrix} \big] (2Z,2\mathbb{B}) \right]_{\epsilon\in\Z^g/2\Z^g},
\end{align*}
\end{definition} 
where $\sigma(Z)=-Z$. This map is an embedding by Lefschetz' theorem \cite{BL}. 

Jacobian varieties are known to be indecomposable, thus we will also require this condition on general $\Rm\in\siegel$, i.e., $\Rm$ must not be symplectically equivalent to a matrix of the form $\mathrm{diag}[\Rm',\Rm'']$ with $\Rm'\in\mathbb{H}^{g'}$, $\Rm''\in\mathbb{H}^{g''}$ and $1\leq g''\leq g'$. Recall that a PPAV $(\Av_\Rm,\Theta_\Rm)$ is indecomposable if there do not exist lower dimensional PPAVs $(\Av_{\Rm'},\Theta_{\Rm'})$, $(\Av_{\Rm''},\Theta_{\Rm''})$ such that $\Av_{\Rm} \cong \Av_{\Rm'}\times \Av_{\Rm''}$ and $\Theta_{\Rm} \cong \Av_{\Rm'}\times\Theta_{\Rm''} + \Theta_{\Rm'}\times\Av_{\Rm''}$ (see \cite{Igu}). 
\begin{definition} \label{def_trisecant}
Let $[V_1]$, $[V_2]$, $[V_3]$ be points in the projective space $\P^n$. They are said to be trisecant or collinear points in $\P^n$ if their representatives $V_1,V_2,V_3 \in \C^{n+1}$ are linearly dependent, i.e., there exist non-zero $c_1,c_2 \in \C$ such that $V_3=c_1 V_1 + c_2 V_2$. 
\end{definition}  
In particular, we say the Kummer variety of $(\Av_\Rm,\Theta_\Rm)$ admits trisecant points if there exist $\Kum(X)$, $\Kum(Y)$, $\Kum(Z)\in\P^{2^g-1}$ satisfying Definition \ref{def_trisecant}.
 
\section{Numerical approaches}
In this section, we outline the numerical approaches to be used in the 
context of the Schottky problem.

\subsection{Computation of theta functions}
The standard way to compute 
the series~(\ref{theta}) is to approximate it by a 
sum, $|N_{i}|\leq \mathcal{N}_{\delta}$, $i=1,\ldots,g$, where the 
constant $\mathcal{N}_{\delta}$ is chosen such that all omitted 
terms in (\ref{theta}) are smaller in modulus than the aimed-at accuracy
$\delta$.  In contrast to~\cite{deconinck03} and \cite{ACh}, we do not give a specific 
bound for each $|N_{i}|$, $i=1,\ldots,g$, i.e., we sum over a 
hyper-cube of dimension $g$ instead of an ellipsoid. The reason for this is that 
it does  not add much to the computational cost, but that it simplifies a 
parallelization of the computation of the theta function in which we 
are interested. Taking into account that we can choose $\mathrm{z}$ in the 
fundamental domain of the Jacobian because of (\ref{eq:periodicity}), 
we get for the Riemann theta function the estimate
\begin{equation}
    \mathcal{N}_{\delta}> \sqrt{-\frac{\ln \delta}{\pi 
    y_{min}}}+\frac{1}{2}
    \label{ne}.
\end{equation}
Here $y_{min}$ is the length of the shortest vector in the lattice 
defined by the imaginary part of the Riemann matrix $\mathbb{B}=X+\mathrm{i} Y$: 
let  $Y=T^\top T$, i.e., $T$ be the Cholesky decomposition of $Y$, then 
$T$ defines a lattice, i.e., a
discrete additive subgroup of $\mathbb{R}^{g}$, of the form 
\begin{equation}
    \mathcal{L}(t_{1},\ldots,t_{g}) = 
    \left\{ TN\bigm| N\in\mathbb{Z}^{g}\right\}
    \label{lattice},
\end{equation}
where $T = [t_1,t_2,\ldots,t_g] \in \mathbb{R}^{g\times g}$ has rank 
$g$. The length of the shortest vector in this lattice is denoted by 
$y_{min}$. 

The greater the norm of the shortest lattice vector, the more 
rapid will be the convergence of the theta series. Note that, in 
general, the convergence of a theta series contrary to popular belief 
can be very slow, see for instance the discussion in 
\cite{deconinck03}. Siegel showed that the length of $y_{min}\geq 
\sqrt{3}/2$ in the Siegel fundamental domain. Unfortunately, no 
algorithm is known to construct a symplectic transformation for a 
general Riemann matrix to this fundamental domain. But Siegel 
\cite{Sie} gave an algorithm to achieve this approximately. A problem 
in this context is the Minkowski ordering. Just as the 
identification of the shortest lattice vector, this is a 
problem for which only algorithms are known whose time grows 
exponentially with the dimension $g$. Therefore, the implementation of 
the Siegel algorithm in \cite{deconinck03} uses an approximation 
to these problems known as the LLL algorithm \cite{LLL}. Whereas this 
algorithm is fast, it is not very efficient. Therefore, in \cite{FJK}, 
Siegel's algorithm was implemented via an exact identification of the 
shortest lattice vector which leads to a more efficient computation 
of the theta functions. Thus, for all Riemann matrices to be 
considered in this paper we always first apply the method outlined in
\cite{FJK} in order to obtain faster convergent theta series. In 
practice, this means that an accuracy of the order $\delta\sim 
10^{-12}$ can be reached with $\mathcal{N}_{\delta}=5$. 

Note that derivatives of theta functions are computed in an analogous 
way as the theta function itself since 
\begin{equation}\label{thetader}
    \frac{\partial\Theta \big[\begin{smallmatrix} \mathrm{p} \\ \mathrm{q} \end{smallmatrix} \big] (\mathrm{z},\mathbb{B})}{\partial \mathrm{z}_{i}}=
    2\pi \I\sum\limits_{\mathrm{N}\in\mathbb{Z}^g}(N_{i}+p_{i})\exp\left\{
    \pi \I \left\langle \mathrm{N}+\mathrm{p}, \mathbb{B}\left(\mathrm{N}+\mathrm{p}\right)
    \right\rangle+2\pi \I
    \left\langle \mathrm{N}+\mathrm{p}, \mathrm{z}+\mathrm{q}
    \right\rangle\right\}
    \;.
\end{equation}
Let $\mathcal{T}(z)$ be the array
\[
  \mathcal{T}(z)=\exp\left\{
    \pi\I \left \langle \mathrm{N}+\mathrm{p}, \mathbb{B}\left(\mathrm{N}+\mathrm{p}\right)
    \right\rangle+2\pi \I
    \left\langle \mathrm{N}+\mathrm{p}, \mathrm{z}+\mathrm{q}
    \right\rangle\right\}, \quad \mathrm{N} \in [-\mathcal{N}_\delta,\mathcal{N}_\delta]^g
\]
The construction of this array is the most expensive part of the routine that computes the theta function, but it allows the immediate computation of the gradient. Let us define the $(2\mathcal{N}_\delta+1)^g \times g$ matrix $M_{\mathrm{N}}$ by
\begin{align*}
M_{\mathrm{N}} = [\mathrm{N}^{(1)}+\mathrm{p}_1, \cdots , \mathrm{N}^{(g)} +\mathrm{p}_g ]
\end{align*}
as the matrix whose row vectors contain all the $\mathrm{N}$ considered in the sum approximating the theta function plus the characteristic $\mathrm{p}$. We use $\mathrm{N}^{(i)}$ to denote the column vector of length $(2\mathcal{N}_{\delta}+1)^{g}$ containing the $i$-th components of all $\mathrm{N}$ in the same order as in the definition of $\mathcal{T}(z)$. Thus, the gradient (given as a row vector) is easily obtained through the matrix multiplication
\begin{equation} \label{grad_theta}
\nabla \Theta \big[\begin{smallmatrix} \mathrm{p} \\ \mathrm{q} \end{smallmatrix} \big](\mathrm{z},\mathbb{B}) \approx 2\pi \I \mathcal{T}(z)^\top* M_{\mathrm{N}}.
\end{equation}

\subsection{Newton iteration}
Newton's iteration can be used to find zeros of locally holomorphic functions of the form $f:\C^M\to\C^N$ if the initial iterate $\mathbf{x}^{(0)}\in\C^M$ is sufficiently close to a zero. With the choice $\mathbf{x}^{(0)}$, we find new iterates via
\begin{equation}
	\mathbf{x}^{(n+1)} = \mathbf{x}^{(n)}- 
	\mbox{Jac}(F(\mathbf{x}^{(n)}))^{-1}F(\mathbf{x}^{(n)}),\quad 
	n=0,1,\ldots;
	\label{newton}
\end{equation}
where $\mbox{Jac}(F(\mathbf{x}^{(n)}))$ denotes the Jacobian matrix of $F$. However, the zero set must be discrete, otherwise $\mbox{Jac}(F(\mathbf{x}^{(n)}))$ becomes singular.\\
The function we are interested in is
\begin{equation} \label{fun_F}
\begin{split}
F: U_\Rm^3 & \to \C^{2^g},\\
(X,Y,Z) & \mapsto c_1(X,Y,Z) \vec{\Theta}(Z) +  c_2(X,Y,Z) \vec{\Theta}(Y)-\vec{\Theta}(X),
\end{split}
\end{equation}
where $\vec{\Theta}(Z) := \left( \Theta \big[\begin{smallmatrix} \epsilon/2 \\ 0 \end{smallmatrix} \big] (2Z,2\mathbb{B}) \right)_{\epsilon\in\Z^g/2\Z^g}$ is the representative of the Kummer point $\Kum(Z)$ in $\C^{2^g}$ and $U_{\Rm}$ is the fundamental domain \eqref{fund} of $\Av_{\Rm}$. Additionally, we need to impose the condition that the zeros must be non-trivial.
\begin{definition}[Trivial zeros]
We say that $(X,Y,Z)$ is a trivial zero if $X\neq \pm Z$, $Y\neq \pm Z$ or $X\neq \pm Y$, since $F(X,Y,Z)=0$ for such cases regardless of the nature of $\Rm\in\siegel$.
\end{definition}
With the aforementioned constraint, $F(X,Y,Z)=0$ if and only if $\Kum(X)$, $\Kum(Y)$, $\Kum(Z)$ are trisecant points, none of which coincide. Therefore, Welters-Krichever's result \cite{Wel,Kri} can be stated in terms of $F$.
\begin{theorem}[Trisecant theorem] \label{welters}
Let $(\Av_{\Rm},\Theta_{\Rm})$ be an indecomposable principally polarized Abelian variety. Then, it is the Jacobian of some Riemann surface $\Rs$ of genus $g$ if and only if the function $F$ has non-trivial zeros.
\end{theorem}
This implies that the zero set of the constrained $F$ is either empty or it is the four-dimensional set given by the parametrization \eqref{XYZ}. However, even if $\Rm$ turns out to be in the Jacobi locus, we still have to determine $X,Y,Z$ without an Abel map, since this is generally unknown. 

Since Theorem \ref{welters} only requires the existence of one zero 
to conclude that $\Rm$ is in the Jacobi locus $\Jg$, we can add four constraints to $F$ in such a way that the zero set of the constrained function is discrete, which will allow us to use Newton's iteration \eqref{newton}. A simple way to do this is by considering its intersection with the zero set of four functions of the form $f(V)=V_j-c$, where $V=(X,Y,Z)\in U_{\Rm}^3$ and $c\in\C$. This is equivalent to fixing four components of the vector $V$ in the iterative process. Having chosen starting vectors satisfying the non-triviality condition, we fix their first components so that their iterates $X^{(n)}$, $Y^{(n)}$ and $Z^{(n)}$ remain different along the whole iteration and, in addition, fix any other component, e.g., $X_2$. Thus, we obtain a function of the form $F|_{W}:W\subset \C^{3g-4} \to \C^{2^g}$, where $W$ is the compact subset 
\begin{equation} \label{W}
W=\{  (X,Y,Z)\in U_\Rm^3\subset\C^{3g} \quad | \quad X_1=X_1^{(0)}, X_2=X_2^{(0)}, Y_1=Y_1^{(0)}, Z_1 =Z_1^{(0)} \}.
\end{equation}
For simplicity, we drop the subscript in $F|_{W}$. Let us denote by $\mathbf{x}$ the remaining $3g-4$ components of $V\in U_{\Rm}^3$, i.e., $\mathbf{x}\in W$. The task is to find a possible zero of the function $F(\mathbf{x})$, i.e., to decide whether it can have values smaller than the specified accuracy $\delta$. If this is the case, then the Riemann matrix is regarded as lying on the Jacobi locus within the given numerical precision.

Thus, after choosing some initial vector $V^{(0)}\in U_\Rm^3$, fixing four components to obtain $W$, setting the initial iterate $\mathbf{x}^{(0)}$, we numerically identify a zero of $F$ by applying the standard 
Newton iteration \eqref{newton}. Notice that $\mbox{Jac}(F(\mathbf{x}^{(n)}))$ is composed of $3g-4$ linearly independent components of the 
full Jacobian matrix of \eqref{fun_F}. Since $F$ is locally holomorphic in the vectors $X$, $Y$, and $Z$, 
the Jacobian can be directly computed via the derivatives 
$\frac{\partial F}{\partial \mathbf{x}_{i}}$, $i=1,\ldots,3g-4$. In 
practice we compute the full Jacobian with respect to all components 
of $X$, $Y$, and $Z$; keep the necessary components and then compute the \emph{Newton step} 
$\mbox{Jac}(F(\mathbf{x}^{(n)}))^{-1}F(\mathbf{x}^{(n)})$ with the 
Matlab command `\textbackslash'. This means that the overdetermined linear system 
$\mbox{Jac}(F(\mathbf{x}^{(n)})) \mathbf{y} = F(\mathbf{x}^{(n)})$ is solved 
for $\mathbf{y}$ in a least squares sense. In our context, this has the advantage that 
all equations in (\ref{Fay3}) are satisfied as required to a certain 
accuracy if the iteration converges.

\begin{remark}\label{rem1}
	If only $3g-4$ equations of $F=0$ are used for the Newton 
	iteration, the latter in general converges to a value of 
	$\mathbf{x}$ for which the remaining components of $F$ do not 
	vanish. Thus, it is important that all $2^{g}$ equations are used 
	in the iteration. 
\end{remark}

We recall that a Newton iteration has quadratic convergence which means that if $\tilde{\mathbf{x}}$ is the wanted zero of $F$, then $\|\mathbf{x}^{(n+1)}-\tilde{\mathbf{x}}\| \propto \|\mathbf{x}^{(n)}-\tilde{\mathbf{x}}\|^{2}$, provided $\mathbf{x}^{(n)}$ is in the so-called \textit{basin of attraction}, which is the subset of $W$ for which the iteration converges.  Loosely speaking the number of correct digits in the iteration doubles in each step of the iteration. However, note that the convergence of the Newton iteration is local and, therefore, depends on the starting point $\mathbf{x}^{(0)}$. If this is not adequately chosen, the iteration may fail to converge.\\

The iteration is stopped either when 
$\|\mathbf{x}^{(n+1)}-\mathbf{x}^{(n)}\|<\delta$ or the residual of 
$F$ is smaller than $\delta$ --- in which case we speak of 
convergence --- or else after 100 iterations, in which case the 
iteration is deemed to not have converged. The latter typically 
indicates an inadequate choice of the initial iterate, and the 
iteration can simply be restarted with another initial vector. When 
the iteration stops, we not only check whether the residual of $F$ in 
(\ref{Fay3}) is below the aimed-at accuracy $\delta$ which would 
indicate that the considered Riemann matrix is in the Jacobi locus. We also compute the \emph{singular value decomposition}\footnote{The singular-value decomposition (SVD) of an $m \times n$-matrix $M $ with complex entries is given by $M = U\Sigma V^{\dagger}$; here $U$ is an $m\times m$ unitary matrix, $V^{\dagger}$ denotes the conjugate transpose of $V$, an $n \times n$ unitary matrix, and the $m \times n$ matrix $\Sigma$ is diagonal (as defined for a rectangular matrix); the non-negative numbers on the diagonal of $\Sigma$ are called the singular values of M.} of the matrix formed by the three vectors $\vec{\Theta}(X)$, $\vec{\Theta}(Y)$, $\vec{\Theta}(Z)$. If the modulus of the smallest singular value of this matrix, denoted by $\Delta$ in the following, is smaller than $\delta$, the vectors are linearly dependent and the Riemann matrix is in the Jacobi locus within numerical accuracy.

\begin{remark}
  It is important that we used in (\ref{c}) ratios of theta functions 
  instead of a form of the Fay identity free of denominators. The reason for this is that odd theta functions have co-dimension one zero sets, and the Newton iteration would converge to all factors in front of $\vec{\Theta}(X)$, $\vec{\Theta}(Y)$, $\vec{\Theta}(Z)$ being zero. The same would happen if we obtain the constants from three of the equations (\ref{Fay3}) in the usual definition of linear dependence of vectors. This is the reason why we also check the linear dependence of the vectors via an SVD. The many zeros of the theta functions in the denominators of (\ref{c}) can also lead to a slow convergence of the Newton iteration for the first iterates. 
\end{remark}


Since the convergence of the Newton iteration depends on the choice 
of the initial iterate $\mathbf{x}^{(0)}$, we must consider a strategy to assure convergence. We have the following choices:
\begin{itemize}
\item[(a)] Take initial vectors of the form with
\begin{equation} \label{x0_half}
V^{(0)}= \frac{\ell}{2} \left(\mathrm{p}^{(1)} + \Rm \mathrm{q}^{(1)},\mathrm{p}^{(2)} + \Rm \mathrm{q}^{(2)},\mathrm{p}^{(3)} + \Rm \mathrm{q}^{(3)} \right),
\end{equation} 
where $0 < \ell < 1$ and $\mathrm{p}^{(i)},\mathrm{q}^{(i)} \in \Z^g/2\Zg$ with $\mathrm{q}^{(i)}\neq 0$. We require $\ell\neq 0,1$ since the Jacobian matrix of the Kummer map is singular at half-period vectors. Thus, this choice will prevent the iterates from being too close to such singular points.
\item[(b)] Take completely random vectors $X^{(0)}$, $Y^{(0)}$, $Z^{(0)}$ in the fundamental domain satisfying the non-triviality conditions.
\item[(c)] Take two random initial vectors and the third one arbitrarily close to one of them, e.g., $V^{(0)}=( X^{(0)}, \ell X^{(0)}, Z^{(0)} )$ with $\ell\approx 1$ or $V^{(0)}= ( X^{(0)}, X^{(0)} + \pmb{\varepsilon}, Z^{(0)} )$ with $\Vert \pmb{\varepsilon} \Vert \approx 0$. We cannot simply set $Y=Z$, $X=Y$ or $X=Z$ because of the non-triviality condition. 
However, from Fay's identity, we know that there exists $Y\neq X$ in the neighbourhood of $X$ making $F$ vanish.
\end{itemize}

For reproducibility purposes, we choose initial vectors according to \eqref{x0_half} with specific values of $\ell$. Thus, for a given Riemann matrix $\Rm_0\in\siegel$, we perform the following steps:

\begin{itemize}
\item[(i)] Given $\Rm_0\in\siegel$, look for a modular transformation 
$R$ such that $\Rm=R \cdot \Rm_0$ is approximately in the Siegel 
fundamental domain (the important point is to ensure that the 
shortest lattice vector has length greater than $\sqrt{3}/2$). This is done with the algorithm \cite{FJK}.
\item[(ii)] Set a value for $\ell$ in (\ref{x0_half}), $\Delta \ell$, $n_{\max}$ and the precision $\delta$. For the examples presented here, we used $\ell=0.10$, $\Delta \ell = 0.10$, $n_{\max}=100$ and $\delta=10^{-10}$.
\item[(iii)] Choose starting vectors $X^{(0)},Y^{(0)},Z^{(0)} \in U_\Rm$ in the form \eqref{x0_half} with $\mathrm{p}^{(1)}=\mathrm{q}^{(1)}=e_{g-2}$, $\mathrm{p}^{(2)}=\mathrm{q}^{(2)}=e_{g-1}$, $\mathrm{p}^{(3)}=\mathrm{q}^{(3)}=e_{g}$ and fix the components $X_1$, $X_2$, $Y_1$, $Z_1$. 
\item[(iv)] Set up the function $F:W\subset \C^{3g-4} \to \C^{2^g}$.
\item[(v)] Perform Newton's iteration \eqref{newton} until $ \| \mathbf{x}^{(n+1)} - \mathbf{x}^{(n)} \| < \delta$ or the maximum number of iterations $n_{\max}$ is reached or $\| F(\mathbf{x}^{(n)}) \| < \delta$. Keep the vectors $X^{(n)}$, $Y^{(n)}$, $Z^{(n)}$ in the fundamental domain at every step. 
\item[(vi)] If $ \Delta^{(N)}:=\min(\mathrm{svd}(\vec{\Theta}(X^{N}), \vec{\Theta}(Y^{N}),\vec{\Theta}(Z^{N}) )) < \delta$, where $N$ is the iteration at which the Newton iteration stopped, we conclude that $\Rm$ is in the Jacobi locus with precision $\delta$ and stop the computations.
\item[(vii)] If the iteration did not converge, replace  $\ell$ by $\ell+\Delta\ell$ and go back to step (iii) and perform a new iteration. Stop if $\ell>\ell_{\max}$ and conclude with a precision $\delta$ that $\Rm$ is not in the Jacobi locus. We used $\ell_{\max}=0.5$ for our examples.
\end{itemize}
We show these steps more explicitly in the pseudo code Algorithm \ref{alg_trisecant}.


\section{Examples in genus 4}
In this section we study known examples in genus 4 since the 
Schottky-Igusa form gives the identification of the Jacobi locus in 
genus 4. This provides interesting tests for our approach. We start 
by studying the stability of the code by considering the Riemann 
matrix of Bring's curve, as well as images of the Abel map (which are 
computed numerically). We use the algorithm described above, but with perturbations of the Abel images as initial iterates. We expect quadratic convergence from the beginning since the initial iterates are in the vicinity of zeros.

\subsection{Bring's curve}
Bring's curve is the curve with the highest number of 
automorphisms in genus 4, see \cite{BN} for the computation of its 
Riemann matrix. It can be defined by the algebraic curve
\begin{equation}
    \{ (x,y)\in\mathbb{C}^{2} \quad | \quad xy^5+x+x^2y^2-x^4y-2y^3=0 \}
    \label{bring}.
\end{equation}
For a given algebraic curve, the Riemann matrix can be computed via 
the symbolic approach by Deconinck and van Hoeij \cite{deco01} implemented in 
Maple or in Sage \cite{SD} or the purely numerical approach 
\cite{FK}, see also the respective chapters in \cite{RSbook}. We use 
here the approach of \cite{FK} and find after applying Siegel's 
algorithm in the form \cite{FJK} the Riemann matrix
\begin{verbatim}
RieMat = 

  -0.5000 + 0.8685i  -0.0000 + 0.0649i  -0.5000 - 0.2678i   0.5000 - 0.2678i
   0.0000 + 0.0649i  -0.5000 + 0.8685i   0.5000 + 0.2678i  -0.5000 + 0.2678i
  -0.5000 - 0.2678i   0.5000 + 0.2678i  -0.0000 + 1.0714i   0.5000 - 0.2678i
   0.5000 - 0.2678i  -0.5000 + 0.2678i   0.5000 - 0.2678i  -0.5000 + 0.8685i.
\end{verbatim}
Note that the accuracy of the computed matrix is estimated to be 
better than $10^{-14}$, but for the ease of readability, we give only 
four digits here. For this matrix, we get for the Schottky-Igusa form 
$|\Sigma|=1.04\times10^{-15}$, i.e., a value of the order of the 
rounding error as expected. This is to ensure that we are indeed testing with a matrix $\Rm$ in the Jacobi locus.

To test the computation of the vectors $X$, $Y$, $Z$ satisfying the Fay identity, we 
consider the points covering the point $z=2.5$ in the complex plane on the sheets 1 to 4 of Bring's 
curve. The code \cite{FK} gives
\begin{verbatim}
AbelMap =

  -0.7052 + 0.3692i   0.0545 + 0.0278i   0.0293 + 0.0775i  -0.0607 + 0.1180i
   0.1286 - 0.2662i   0.2747 + 0.2456i  -0.4068 + 0.4113i   0.0318 + 0.2067i
  -0.4351 + 0.2906i  -0.2108 + 0.0422i   0.0250 + 0.2906i   0.0823 - 0.2451i
   0.4519 - 0.6915i   0.0718 + 0.0915i  -0.0126 - 0.0140i   0.0487 + 0.0493i.
\end{verbatim}
We use these vectors (the columns of the matrix) and set $X^{(0)}$, $Y^{(0)}$, $Z^{(0)}$ as the linear combinations given by \eqref{XYZ} for this particular example, rather than using the choice \eqref{x0_half} as we will do in general.
The residual of the function $F$ in (\ref{Fay3}) for this Abel map is of 
the order of $10^{-10}$, and the minimal singular value of the matrix 
with the vectors $\Vec{\Theta}(X)$, $\Vec{\Theta}(Y)$ and $\Vec{\Theta}(Z)$ is $\Delta\sim 
10^{-11}$. This indicates that the Abel map is computed to an 
accuracy of the order of $10^{-11}$. After one Newton iteration, the 
difference between the new and old vector is of the order $10^{-11}$, 
and the residual of $F\sim 10^{-14}$ and $\Delta\sim10^{-14}$.  This 
shows that the iteration is stable (the numerical error in the Abel 
map is `corrected' by the Newton iteration), and that a similar 
residual of $F$ is reached in this example as for the Schottky-Igusa 
form.

The stability of the iteration is also shown by perturbing the initial vector. We multiply the 
above $\mathbf{x}$ by a factor $1.1$, i.e., we keep the entries $X^{(0)}_1$, $X^{(0)}_2$, $Y^{(0)}_1$, $Z^{(0)}_1$ and multiply the remaining ones by $1.1$. After 7 iterations, we get 
the same residual, minimal singular value $\Delta$ and final vector $\mathbf{x}$ as before. The same behaviour is observed for $0.9 \mathbf{x}$ as the initial iterate.

It is known that 
Matlab timings are not very precise since they strongly depend on how 
many precompiled commands are used in the coding, but they provide an 
indication of actual computing times for a given task. On a standard 
laptop, the above examples take a few seconds. 

To finish the tests with Bring's curve, we choose the same $X^{(0)}_1$, $X^{(0)}_2$, $Y^{(0)}_1$, $Z^{(0)}_1$ as before, but the remaining components are chosen randomly with the condition that $X^{(0)}$, $Y^{(0)}$, $Z^{(0)}$ must be in the fundamental domain. The algorithm finally produces a residual for $F$ smaller than $10^{-12}$, but finds a vector $\mathbf{x}$ different from the one produced by the Abel map. This is because the zero set of $F|_W$, with $W$ given by \eqref{W}, does not necessarily have a unique element.

\subsection{Family of genus 4 Riemann matrices}
We now turn to the family of Riemann matrices studied in 
\cite{CKS}, which correspond to the genus 4 algebraic curves
\begin{equation} \label{GM_curve}
    \{ (x,y)\in \C^2 \quad | \quad y^2=x(x+1)(x-t) \},
\end{equation}
parametrized by $t\in \P^1 \backslash \{0,1,\infty \}$. Their period matrices are given by \cite{GM}, i.e., $\Rm_{\tau} = A^{-1} B$, where
\begin{equation}
    (A | B)=\begin{pmatrix} \tau & \tau& 0& -\tau-1& 1& 1& 0& -1 \\
    \zeta^2-1 & 1& -\zeta^2+1& 1& 1& -\zeta^2& \zeta^2& -\zeta^2+1 \\
    \zeta^3-\zeta& -\zeta^3& -2\zeta^3+2\zeta^2+\zeta-1& \zeta^2-\zeta& 1& \zeta^2-1& \zeta^3-\zeta^2-2\zeta+2& \zeta^2 \\
    -\zeta^3+\zeta &\zeta^3& 2\zeta^3+2\zeta^2-\zeta-1& \zeta^2+\zeta& 1& \zeta^2-1& -\zeta^3-\zeta^2+2\zeta+2& \zeta^2
    \end{pmatrix},
\end{equation}
with $\zeta=e^{2\pi\I/12}$ and $t=\mu(\tau)$ for some $\mu: \mathbb{H}\to \P^1$.

\begin{figure} [H]
\includegraphics[width=6cm]{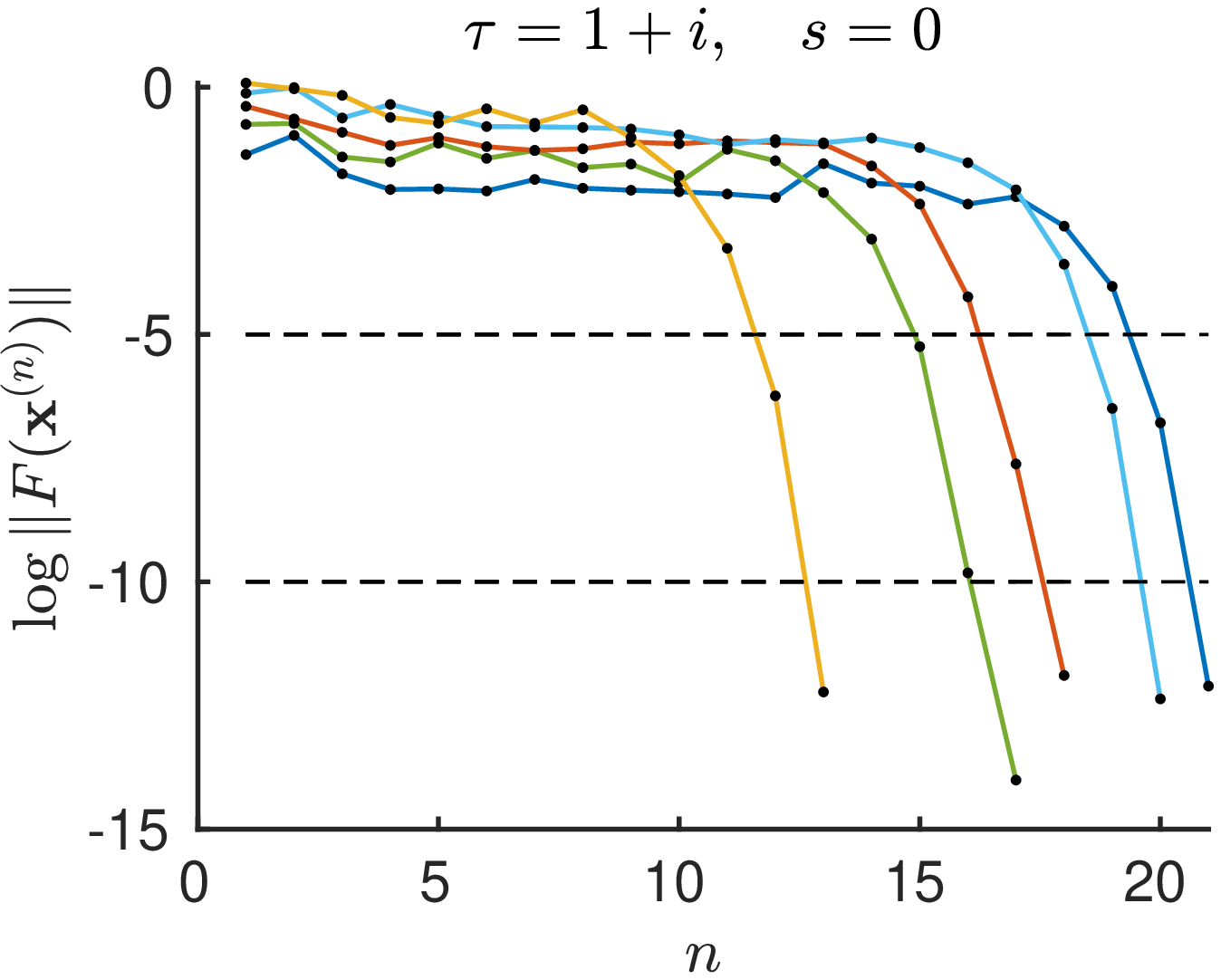} \includegraphics[width=6cm]{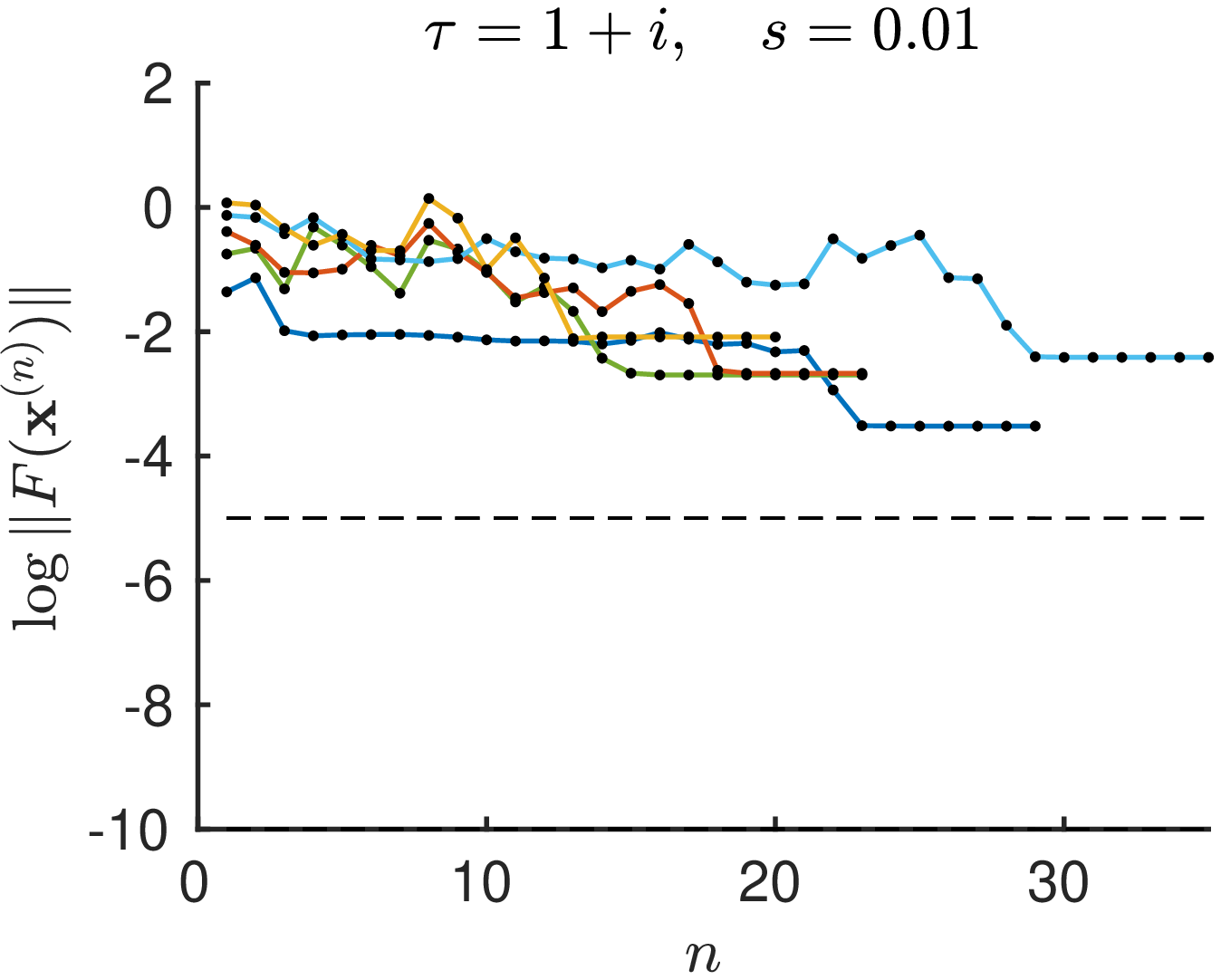}
\caption{Residual of the function $F$, corresponding to the matrix $\Rm_{\tau,s}$, in dependence of the $n$-th iteration.}  \label{fig_conv_genus4}
\end{figure}

\begin{figure} [H]
\includegraphics[width=6cm]{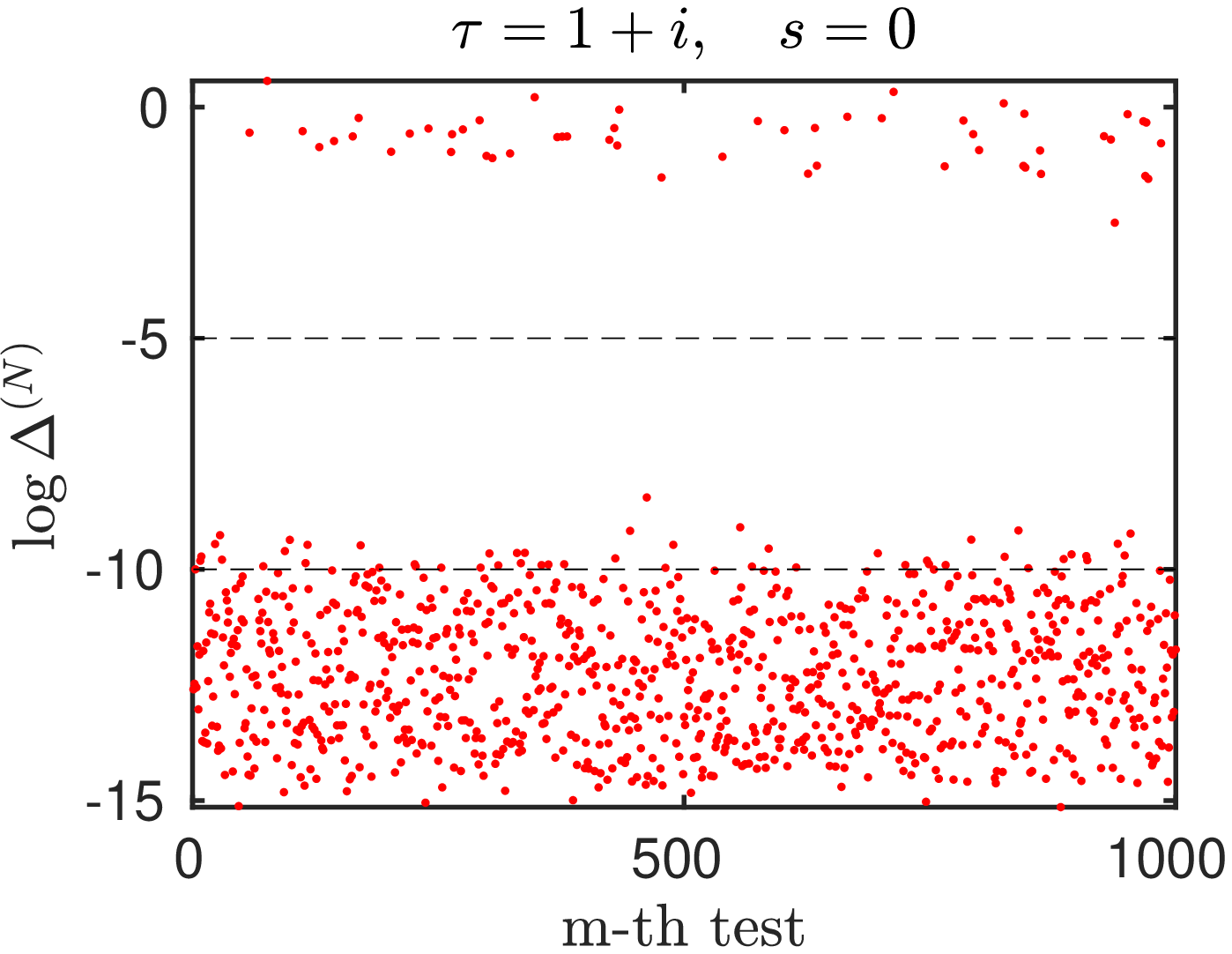} \includegraphics[width=6cm]{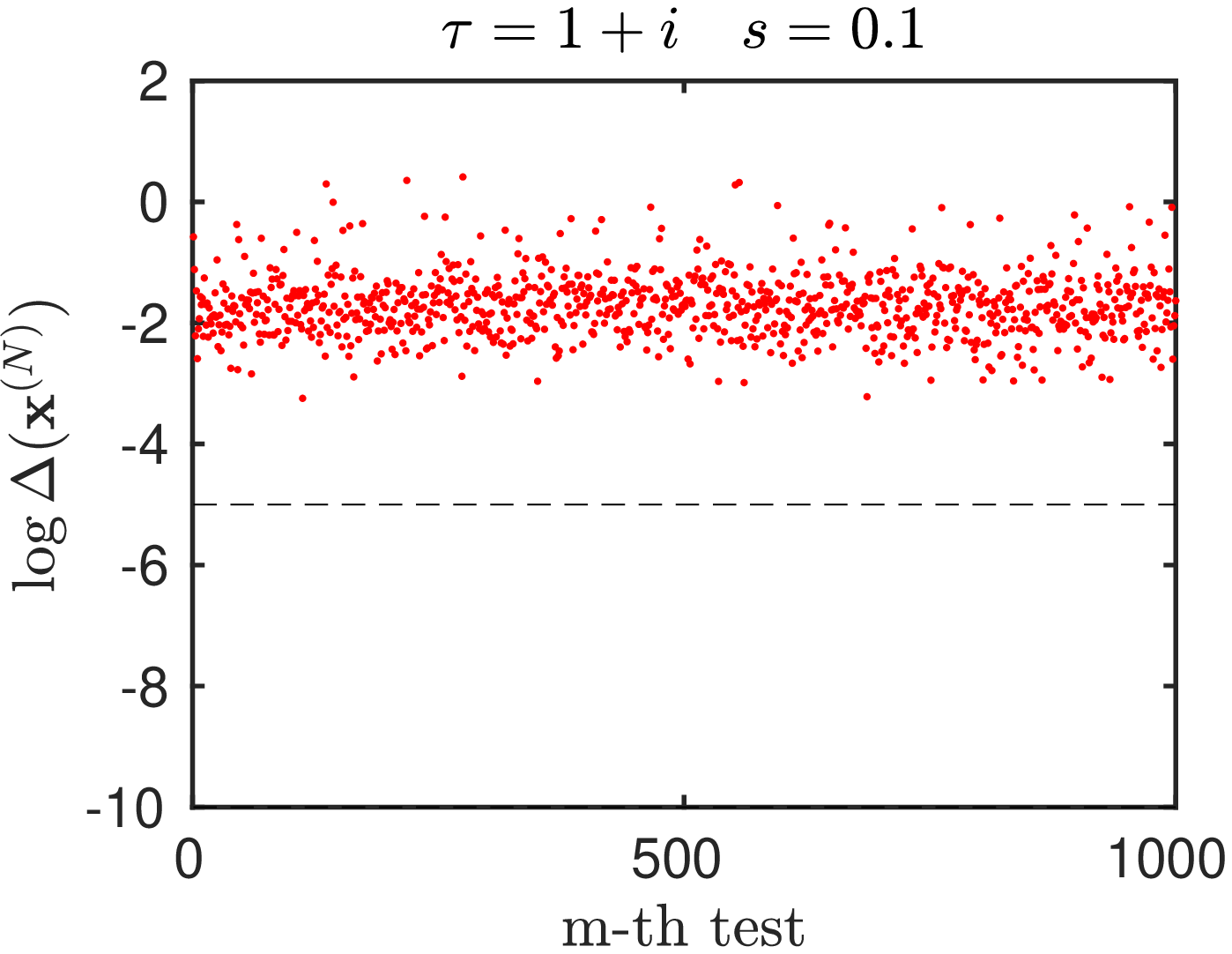}
\caption{Attained values of $\Delta^{(N)}$ for randomly chosen initial vectors. The left-hand side corresponds to $\Rm_{\tau,s}$ with $s=0$ and the right-hand side to $s=0.1$.}  \label{fig_residuals_genus4}
\end{figure}

Let us first look at the convergence rate of the iteration for one specific Riemann matrix, e.g., $\Rm_\tau$ with $\tau=1+\I$. In Figure \ref{fig_conv_genus4} we observe that as soon as the vector $\mathbf{x}^{(n)}$ falls into a basin of attraction, it converges to a zero very rapidly. Recall that $N$ is the step at which the iteration stops. Thus, $\Vert F(\mathbf{x}^{(N)})\Vert $ is the best residual achieved by the iteration corresponding to the starting vectors $X^{(0)}, Y^{(0)}, Z^{(0)}$. For this test, we chose the starting vectors according to Algorithm \ref{alg_trisecant}. In contrast, notice that with a small perturbation of the form 
\begin{equation} \label{B_s}
\Rm_{\tau,s}=\Rm_\tau + s \cdot \mathrm{diag} [2,3,5,7], \quad \text{with } s=0.01;
\end{equation}
the smallest value of $\Vert F(\mathbf{x}^{(N)})\Vert $ is above $10^{-4}$ (although this lower bound will depend on how small $s$ is). Thus, for these particular examples there are at least six orders of magnitude of difference in the smallest attained residuals.
\begin{remark}
Using Algorithm \ref{alg_trisecant}, but with randomly chosen initial vectors, the iterations converge with a frequency of approximately $90 \%$ when $\B$ is in the Jacobi locus. We attain the residuals shown in Figure \ref{fig_residuals_genus4} for 1000 different tests for the same matrix $\Rm_\tau$, i.e., 1000 iterations with different randomly chosen initial vectors $X^{(0)}, Y^{(0)}, Z^{(0)}$ attained $\Delta^{(N)} < 10^{-10}$ 896 times. In contrast, none of the tests attained a value below $10^{-5}$ when the considered matrix is  $\Rm_{\tau,s}$ with $s=0.1$.
\end{remark}

\begin{figure} [H]
\includegraphics[width=6cm]{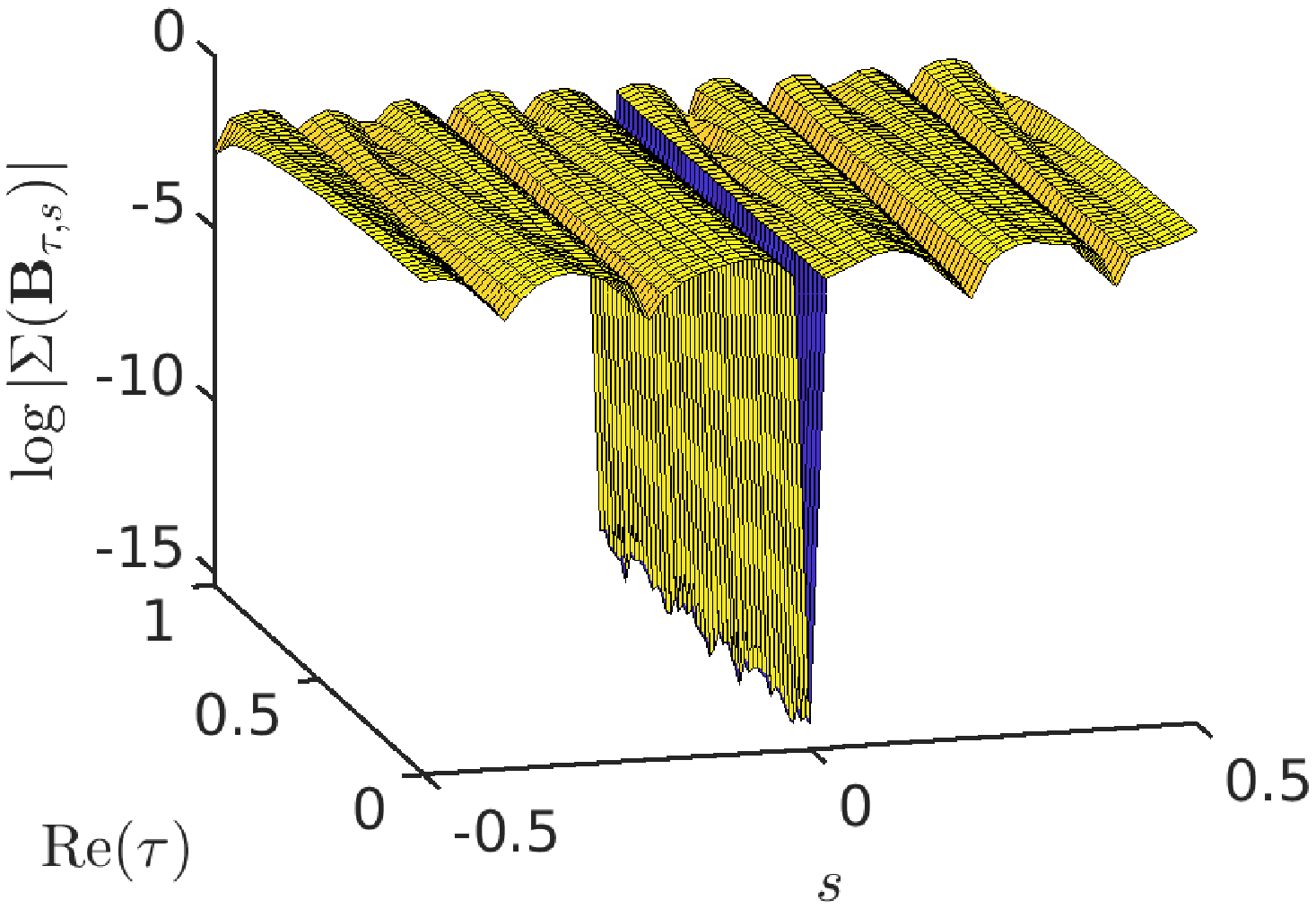} \includegraphics[width=6cm]{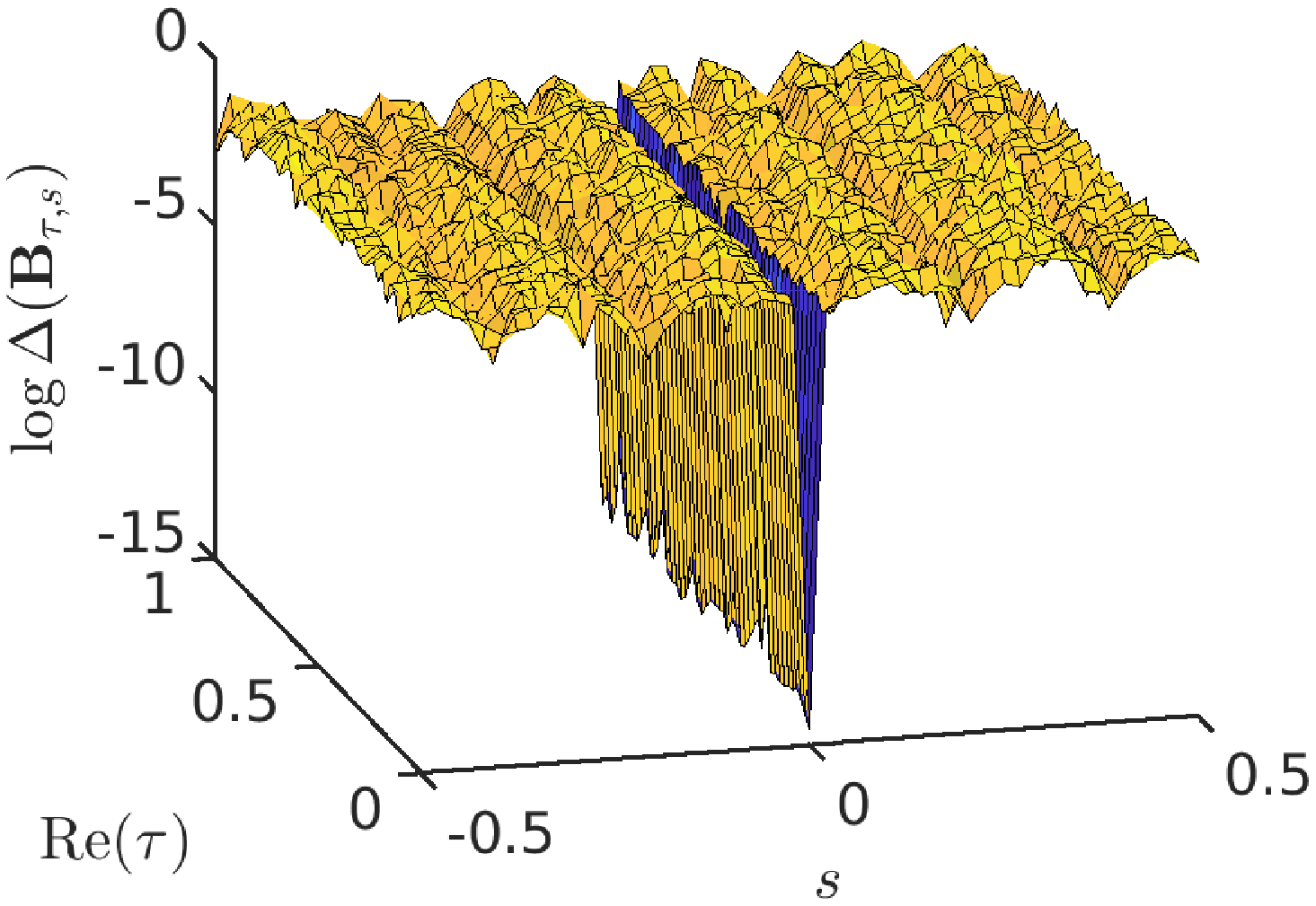} 
\caption{On the left the Schottky-Igusa form and on  the right the 
minimum singular value $\Delta^{(N)}$, both in dependence of the Riemann matrices $\Rm_{\tau,s}$.}  \label{fig_SI_form}
\end{figure}

Let us consider the family of Riemann matrices of the form \eqref{B_s} parametrized by $s\in[-1/2,1/2]$ and $\tau=x+\mathrm{i}$, with $x\in[0,1]$. We expect the Schottky-Igusa form of every $\Rm_{\tau,s}$ and the residual $\| F(\mathbf{x}^{(N)}) \|$ (or its associated minimum singular value) obtained with Algorithm \ref{alg_trisecant} to be within the same order of magnitude. This is indeed what we observe in Figure \ref{fig_SI_form}. With this algorithm, we conclude with precision $\delta=10^{-10}$ that the Riemann matrices $\Rm_{\tau,s}$ with $s\neq 0$ are not in the Jacobi locus, agreeing with the Schottky-Igusa form.

Note that our approach produces numerically essentially the 
same residuals as the Schottky-Igusa form and is thus of the same 
practical relevance in deciding whether a matrix is in the Jacobi 
locus or not with given numerical accuracy. In both cases the input 
is only the Riemann matrix $\mathbb{B}$, but in the Fay approach an initial 
iterate for the vectors corresponding to Abel maps for $\mathbb{B}$ 
in the Jacobi locus.

The matrices $\Rm_{\tau}$ considered above are in exact form, thus the smallest residual is expected to be zero to machine precision, i.e., approximately $10^{-14}$ with Matlab. However, the smallest $\Vert F(\mathbf{x}^{(N)}) \Vert$ might be considered as an indicative of the precision of its input matrix $\Rm$, if this is not in exact form. For example, let us  consider the matrices $ \Rm_{\tau,s} = \Rm_\tau +s (M+\mathrm{i} M) $ where $M$ is the symmetric matrix with coefficients $M_{jk}=(j+k)/5$, $\tau=1+\mathrm{i}$ and $s\in [10^{-15},10^{-1}]$ are small perturbations. Thus, $\Rm_{\tau,s}$ can be seen as the Riemann matrices of the curve \eqref{GM_curve} with an accuracy of the order $s$. We remove the stopping criterion $\Vert F(\mathbf{x}^{(N)}) \Vert < \delta$ for this particular example in order to visualize the smallest residual we can achieve with Matlab's precision.
\begin{figure} [H] \label{logS}
\includegraphics[width=7cm]{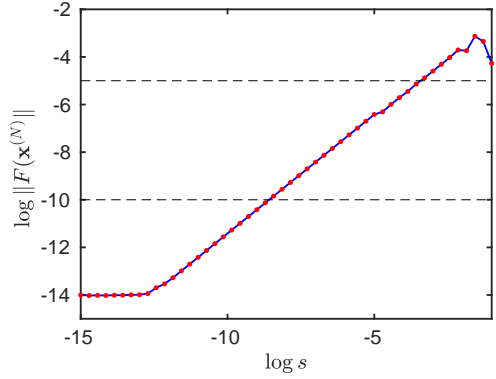}
\caption{Smallest residuals of period matrices with accuracy $s$.}
\end{figure}
We can observe in the figure that we attain a residual to machine precision if the matrices are within a precision of approximately $10^{-13}$ and, in particular, the residual is below $\delta=10^{-10}$ if the Riemann matrix is within a precision of the order $\delta$. This is important, since some of the Riemann matrices in higher dimensions are only known numerically.

\section{Examples in higher genus}
We perform similar tests with matrices in higher genus such as, for example, the period matrices of hyperelliptic curves of arbitrary genus given by \cite{Schi}. In this paper we only show examples corresponding to the curve
\begin{align} \label{hyper_curve}
\{ (x,y)\in\C^2 \quad | \quad  y^2 = x(x^{2g+1}-1)  \},
\end{align}
but similar results are obtained with the period matrices of other 
hyperelliptic curves. The period matrix of \eqref{hyper_curve} is  
\begin{equation}
\B_{jk} = 1+\frac{1}{\tau_1} \sum_{l=1}^j \tau_l \tau_{k-j+l} \quad \text{for } 1\leq j \leq k \leq g,
\end{equation}
where $\tau_1 = (-1)^g \zeta_{2g+1}^{g^2}$ with  $\zeta_m^n=\exp(2\pi \I m/n) $ and
\begin{align*}
\tau_{j+1} = \frac{\tau_1}{1+\zeta^{j}_{2g+1}} \left( 1 - \sum_{l=2}^j\zeta_{2g+1}^{g-j+l-1} \tau_l \tau_{j-l+2} \right), \quad j=1,...,g-1.
\end{align*}
We use these matrices for tests up to $g=7$, but we also consider period matrices of non-hyperelliptic curves. For $g=6$  the Fermat curve 
\begin{equation} \label{fermat}
\{ [X:Y:Z]\in \mathbb{P}^2 \quad | \quad X^m+Y^m=Z^m \},
\end{equation}
with $m=5$, provides an example. With the numerical 
approach \cite{FK}, we get after applying the algorithm \cite{FJK} 
the Riemann matrix
\begin{verbatim}
	RieMat =

  Columns 1 through 4

  -0.3735 + 0.9276i  -0.3574 + 0.4580i  -0.4578 + 0.3092i  -0.2891 + 0.3705i
  -0.3574 + 0.4580i   0.1365 + 1.0006i  -0.0161 + 0.4697i   0.1104 - 0.1415i
  -0.4578 + 0.3092i  -0.0161 + 0.4697i   0.3474 + 1.0079i  -0.2630 - 0.3894i
  -0.2891 + 0.3705i   0.1104 - 0.1415i  -0.2630 - 0.3894i  -0.3152 + 1.1305i
   0.0905 + 0.4390i  -0.4616 + 0.4201i   0.3635 + 0.5382i  -0.3735 - 0.2479i
  -0.4417 - 0.1605i  -0.3313 - 0.3020i   0.2891 - 0.3705i  -0.1725 + 0.0496i

  Columns 5 through 6

   0.0905 + 0.4390i  -0.4417 - 0.1605i
  -0.4616 + 0.4201i  -0.3313 - 0.3020i
   0.3635 + 0.5382i   0.2891 - 0.3705i
  -0.3735 - 0.2479i  -0.1725 + 0.0496i
  -0.4839 + 1.0692i  -0.3796 - 0.0685i
  -0.3796 - 0.0685i  -0.4095 + 0.8023i.
\end{verbatim}
For this curve, we could again compute an Abel map, but we are 
interested also in perturbations of this Riemann matrix not in the 
Jacobi locus. \\
Analogously to \eqref{B_s}, we add diagonal perturbations of the form $\Rm_s = \B + s\cdot \mathrm{diag}[2,...,g+1]$ with $s\in[-1/2,1/2]$.  As in the previous example, we observe several orders of magnitude of difference between the case $s=0$ and the cases $s\neq 0$. Although a similar behaviour to Figure \eqref{logS} is to be expected for small values of $s$.
\begin{figure} [H]
\includegraphics[width=6cm]{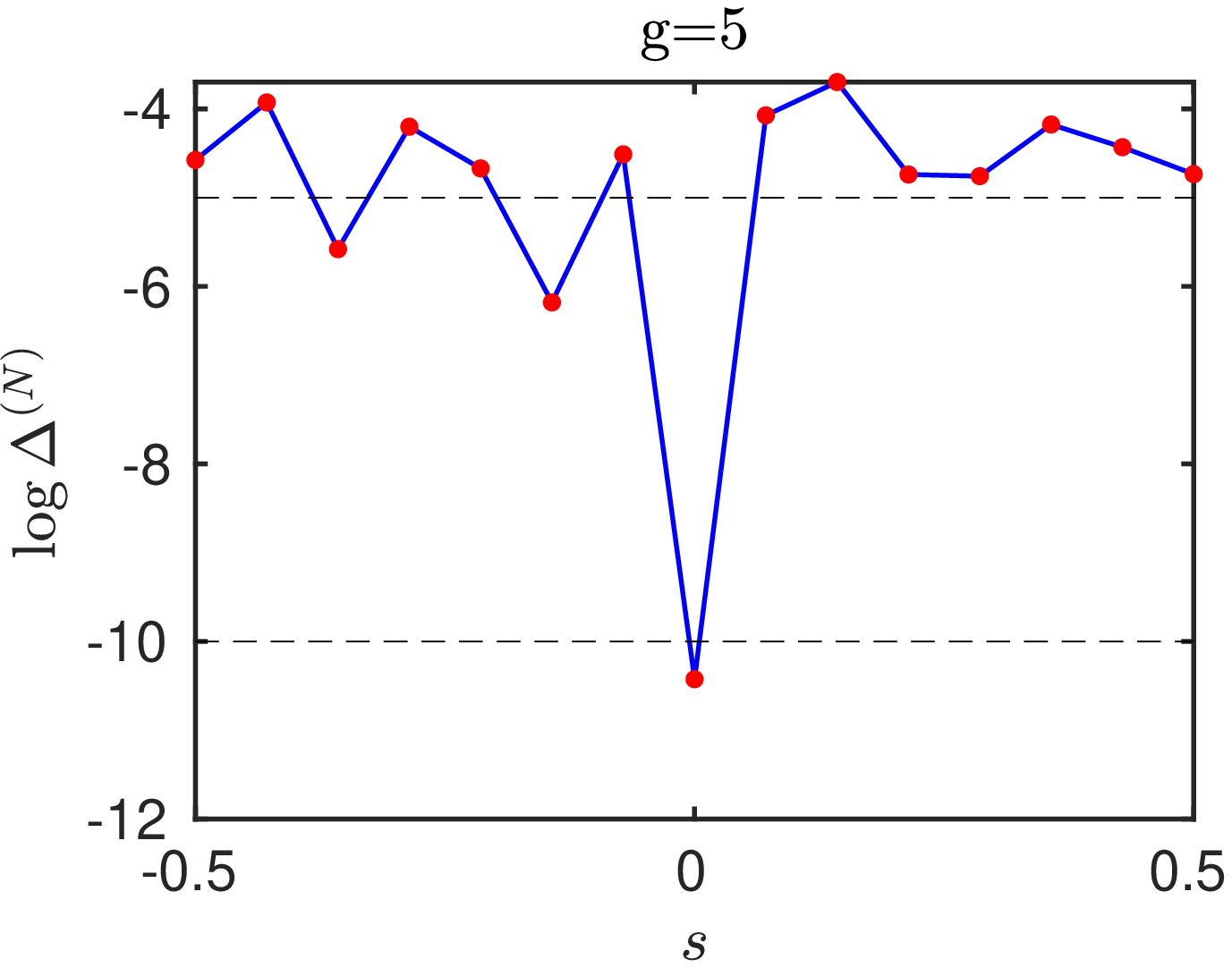} \includegraphics[width=6cm]{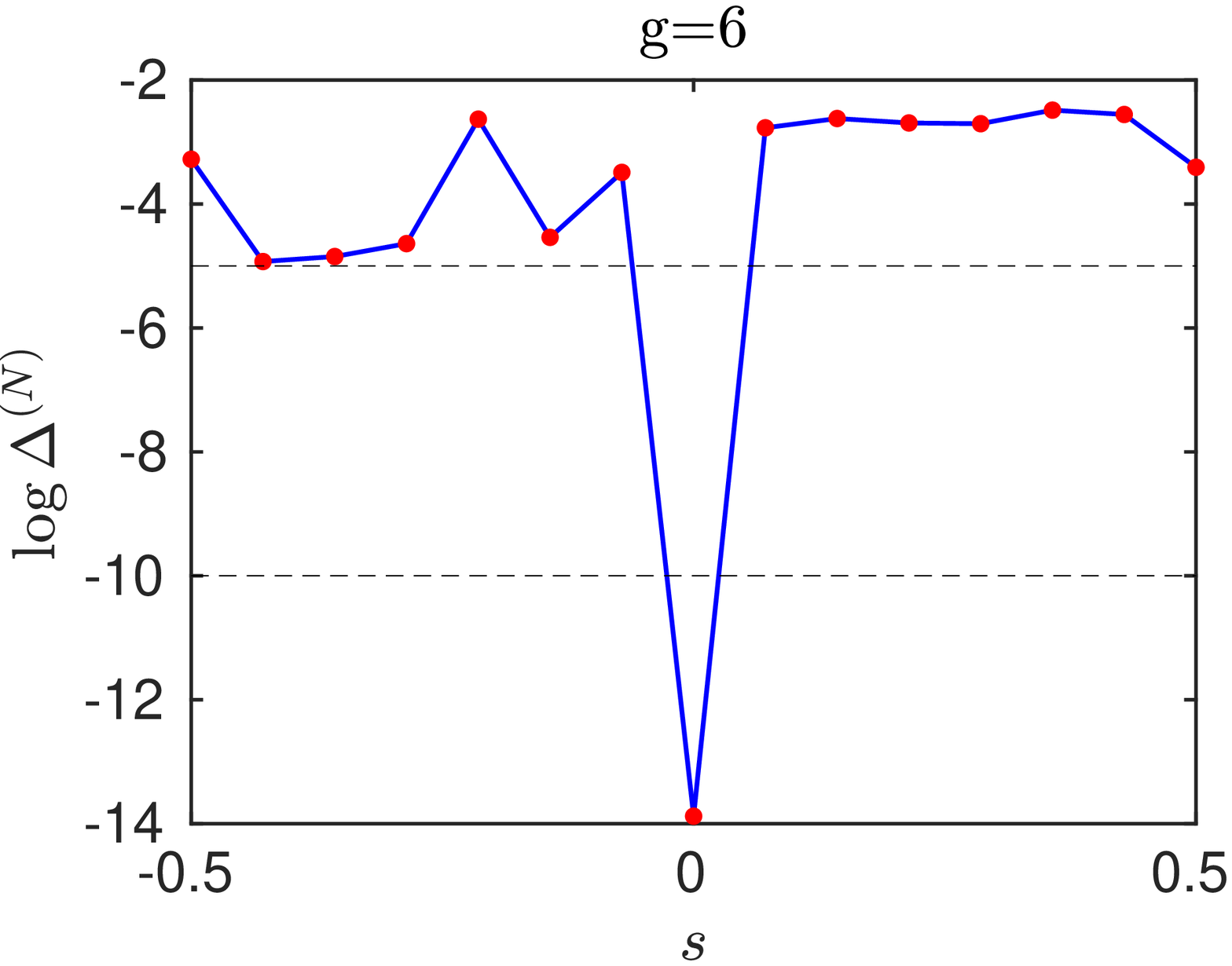}
\caption{Minimum $\Delta ^{(N)}$ obtained for the matrices $\Rm_s\in\siegel$ for $g=5$ and $g=6$.} \label{fig_hyper_g4_g5}
\end{figure}

The Fricke-Macbeath curve \cite{fricke,macbeath} is a curve of
genus 7 with the maximal number $84(g-1)=504$ of
automorphisms. It can be defined via the algebraic curve
\begin{equation}
    \{ (x,y)\in\C^2 \quad | \quad 1+7yx+21y^2x^2+35x^3y^3+28x^4y^4+2x^7+2y^7=0\}.
    \label{FM}
\end{equation}
After applying the algorithm from~\cite{FJK}, the code \cite{FK} leads to 
the Riemann matrix
\begin{verbatim}
RieMat =

  Columns 1 through 4

   0.3967 + 1.0211i   0.0615 - 0.1322i   0.0000 - 0.0000i   0.4609 + 0.2609i
   0.0615 - 0.1322i   0.3967 + 1.0211i  -0.3553 + 0.5828i   0.3386 - 0.1933i
   0.0000 - 0.0000i  -0.3553 + 0.5828i   0.2894 + 1.1656i   0.0905 + 0.2450i
   0.4609 + 0.2609i   0.3386 - 0.1933i   0.0905 + 0.2450i   0.3967 + 1.0211i
  -0.3553 + 0.5828i  -0.4776 + 0.1287i  -0.4776 + 0.1287i  -0.4776 + 0.1287i
  -0.1838 - 0.3219i  -0.2743 - 0.5669i   0.3871 - 0.3736i   0.0167 - 0.3895i
   0.3386 - 0.1933i   0.3386 - 0.1933i  -0.1223 - 0.4541i   0.0615 - 0.1322i

  Columns 5 through 7

  -0.3553 + 0.5828i  -0.1838 - 0.3219i   0.3386 - 0.1933i
  -0.4776 + 0.1287i  -0.2743 - 0.5669i   0.3386 - 0.1933i
  -0.4776 + 0.1287i   0.3871 - 0.3736i  -0.1223 - 0.4541i
  -0.4776 + 0.1287i   0.0167 - 0.3895i   0.0615 - 0.1322i
   0.2894 + 1.1656i  -0.1671 - 0.7115i   0.0905 + 0.2450i
  -0.1671 - 0.7115i   0.4414 + 1.2784i  -0.3386 + 0.1933i
   0.0905 + 0.2450i  -0.3386 + 0.1933i   0.3967 + 1.0211i.
\end{verbatim}
The accuracy of this matrix is estimated to be of the order of 
$10^{-10}$.

\begin{figure} [H]
\includegraphics[width=6cm]{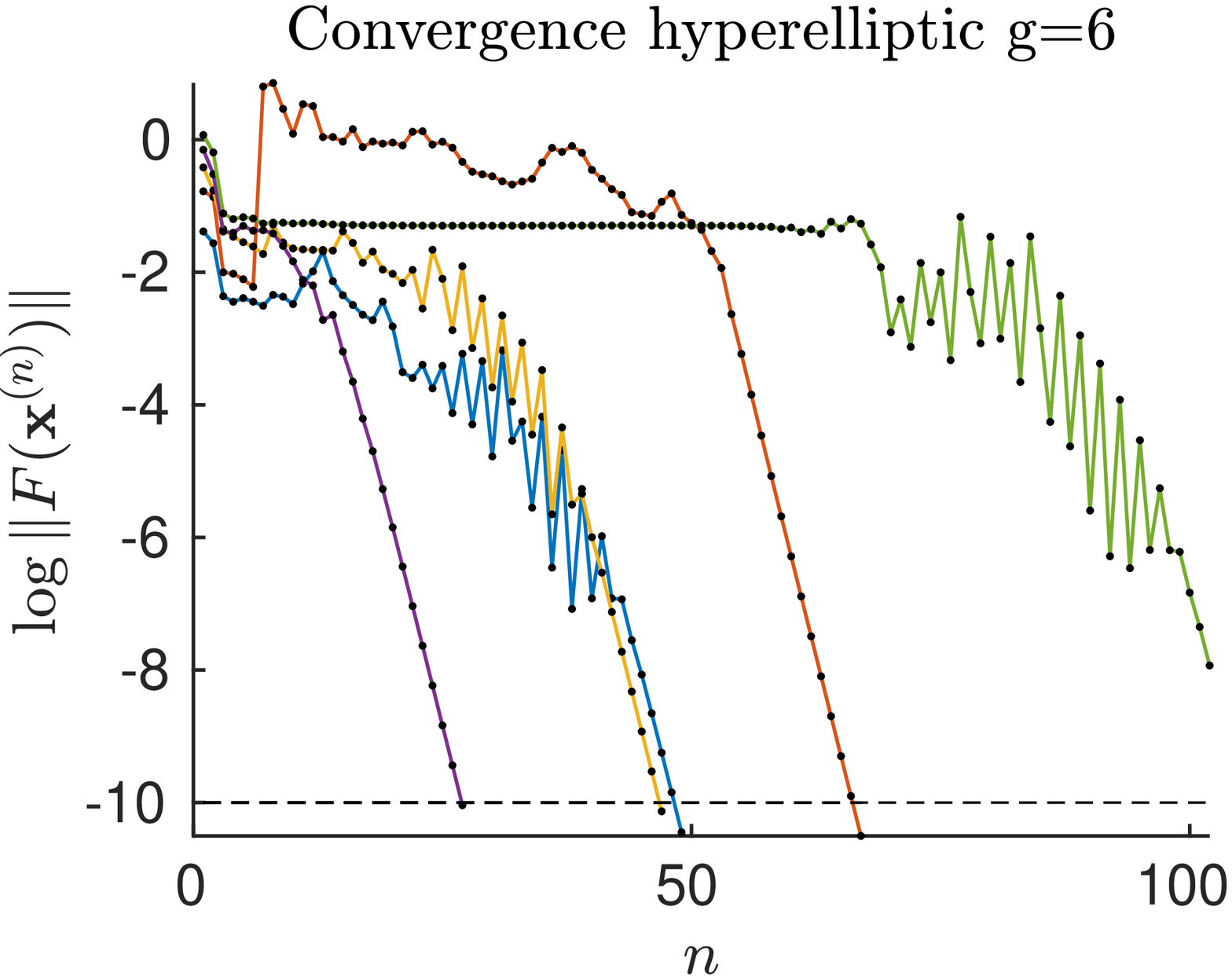} \includegraphics[width=6cm]{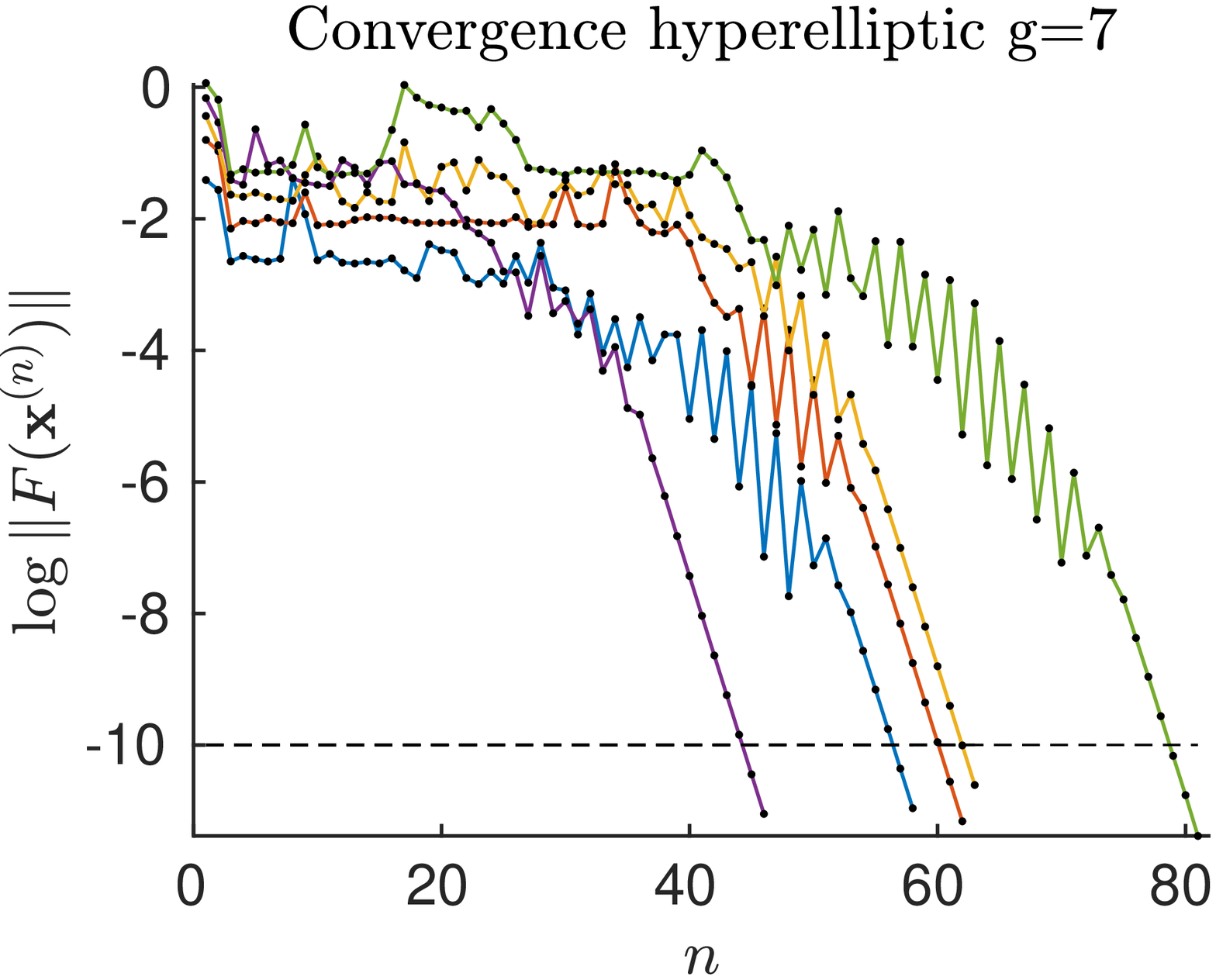}
\caption{Residual of the function $F$ in dependence of the $n$-th iteration.}  \label{fig_hyper_g6_g7}
\end{figure}

Fast convergence when $\mathbf{x}^{(n)}$ falls into a basin of attraction is still observed in higher genus. 
For example, let us observe the convergence of the residual of $F$ corresponding to the period matrix of Fermat curve for $g=6$ and the period matrix of the Fricke curve for $g=7$. The dimension of the domain of $F$ increases linearly with $g$, thus it might take more steps for the iteration to find a basin of attraction, but when it does, fast convergence is assured.
\begin{figure} [H]
\includegraphics[width=6cm]{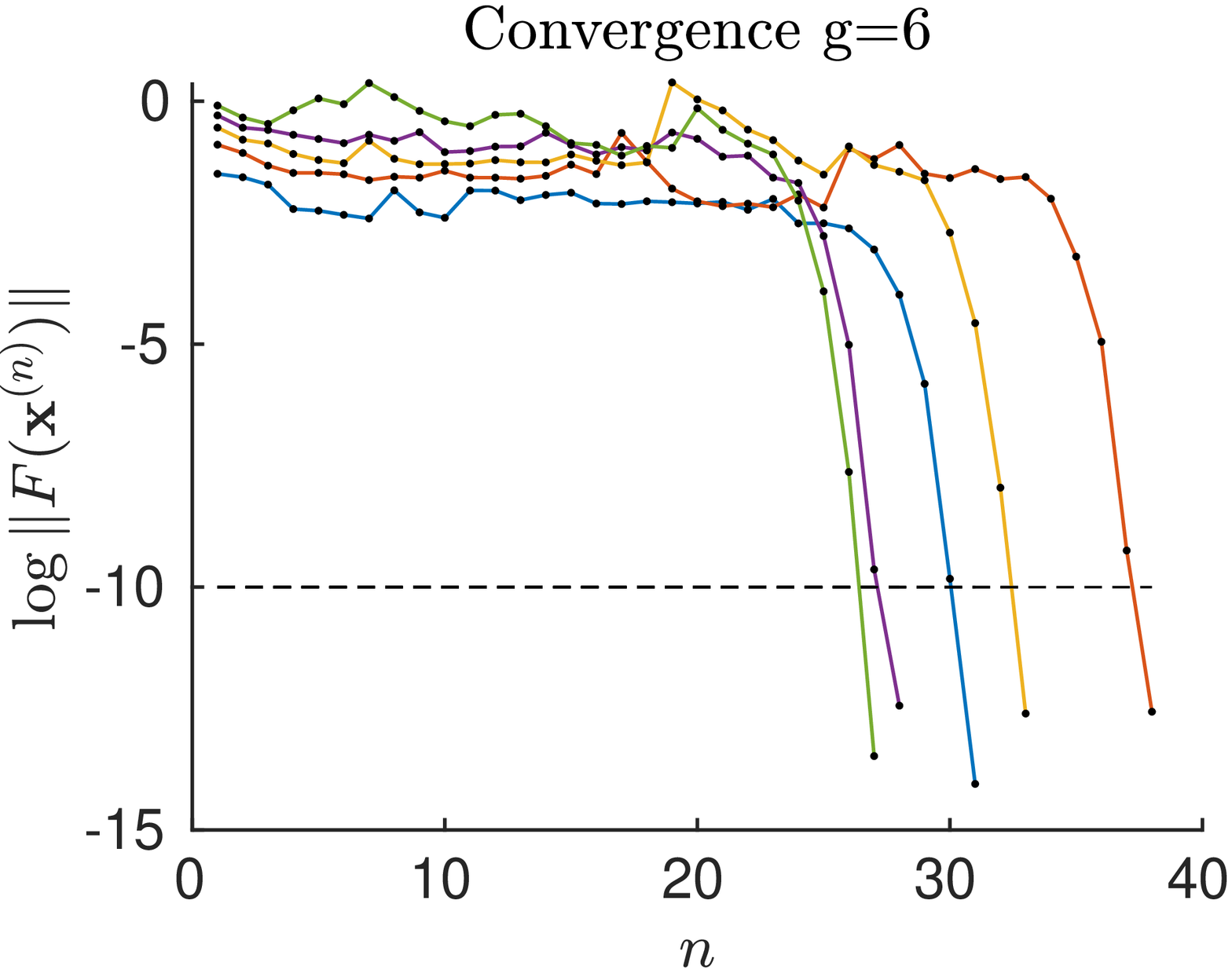} \includegraphics[width=6cm]{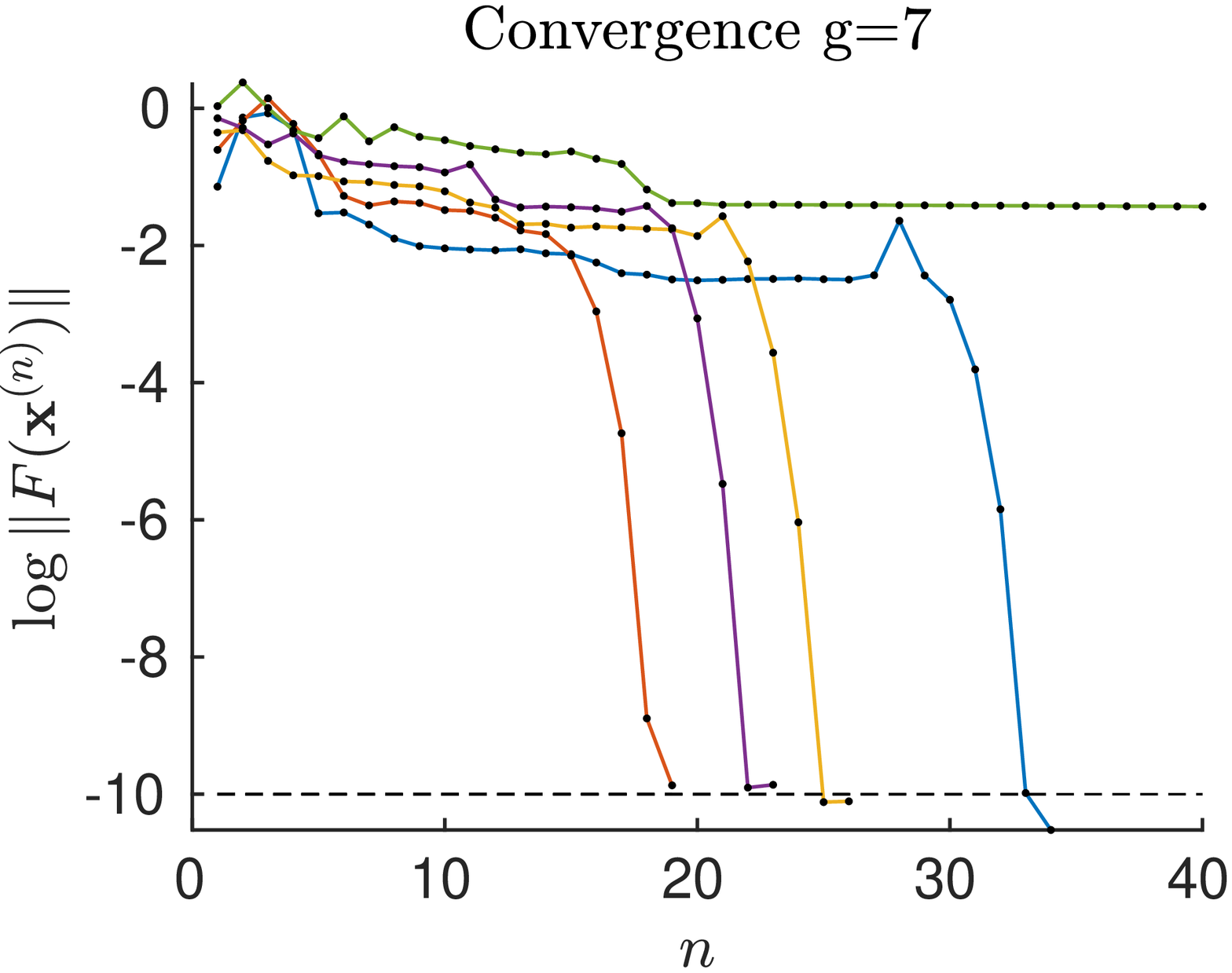}
\caption{Residual of the function $F$ in dependence of the $n$-th iteration.} \label{fig_g6_g7}
\end{figure}
It is worth mentioning that the convergence is slower if the matrix 
$\Rm$ corresponds to the hyperelliptic example, as we can see in Figure \ref{fig_hyper_g6_g7} for the curve \eqref{hyper_curve} with $g=6$ and $g=7$. This might be due to the extra symmetries in their matrix elements, which makes the convergence slower. However, these are just special cases amongst all the possible matrices in $\siegel$.

The computation per step becomes more expensive as $g$ increases since the number of summands in the approximation of every theta function is $(2\mathcal{N}_\delta+1)^g$. Besides that, we need to compute $2^g$ even and $6$ odd theta functions in each step. This increases the computational time of the whole algorithm. As a reference, the total time of the tests for the period matrices of non-hyperelliptic curves with $g=4,5,6$ presented in Figures \ref{fig_conv_genus4} and \ref{fig_g6_g7} consisting of $5$ iterative processes are
\begin{center}
\begin{tabular}{ |c|c| } 
 \hline
 $g$ & time  \\ \hline
 4 & 6 s \\
6 & 28 min \\ 
 7 & 19 h  \\ 
 \hline
\end{tabular}.
\end{center}
However, since it suffices finding only one zero, the computations can also be terminated as soon as we find a vector with $\Vert F(\mathbf{x}^{(n)}) \Vert < \delta$, thus reducing these times. On the other hand, the total computational time of the five iterative processes for the period matrices of hyperelliptic curves are
\begin{center}
\begin{tabular}{ |c|c| } 
 \hline
 $g$ & time  \\ \hline
4 & 14 s \\ 
 5 & 98 s \\ 
 6 & 56 min \\
 7 & 30 h \\
 \hline
\end{tabular}
\end{center}

\section{Outlook}
In this paper, we have shown how the Fay identity can be used to 
identify the Jacobi locus for small genera in a computational 
context.  Concretely, the Fay identity in the form (\ref{Fay3}) is 
considered for a given Riemann matrix $\mathbb{B}$ and its fulfilment 
is probed by using a Newton iteration. Four components of the vectors 
$X$, $Y$ and $Z$ are fixed in order to avoid 
convergence to a trivial solution where two or more of the vectors 
coincide. This leads to an iteration to find the zeros of a function $F$ depending on the vector $\mathbf{x}$ with $3g-4$ complex components.
    The starting point for the 
iteration is the choice (\ref{x0_half}).  

The approach uses a finite precision $\delta$, and the answer to whether a matrix 
is in the Jacobi locus is thus given with this accuracy. In this 
paper, we have only studied the case of so-called double precision 
leading in practice to a $\delta\sim10^{-12}$. However, the code is 
set up in a way that inclusion of multi-precision packages is 
straightforward allowing essentially arbitrary precision. 

We have considered in this paper examples up to $g=7$. This limit is 
imposed by the size of memory on the used standard computers for the computation 
of the theta functions. Since 
the algorithm \cite{FJK} allows the computation of theta functions 
and their derivatives with a truncation parameter 
$\mathcal{N}_{\delta}=5$ in double precision, these computations are very efficient. In order to take advantage of Matlab's 
vectorization algorithms, the exponentials of the theta function 
(\ref{theta}) are stored in one array with $11^{g}$ components. If 
one wants to go to higher values of $g$, the memory limitations can 
be avoided by not storing these values and applying loops to compute 
them.  The Jacobian matrix in (\ref{newton}) is of complex dimension $2^{g}\times (3g-4)$, 
but its 
computation is also fully parallelizable.  It will be the subject 
of further research which genera can be reasonably studied with this 
approach. 

We note that the algorithm produces as a by-product three vectors 
being linear combinations of the Abel maps of 
four points if the Riemann matrix is identified to be in the Jacobi 
locus. The resulting vectors are just auxiliary to test whether the 
Fay identity is satisfied. It is an interesting question whether they 
can be used to identify the defining algebraic curves in a similar way as 
in \cite{ACE}. 

\appendix
\section{Algorithm}
\begin{algorithm}[H]
\caption{Algorithm to find trisecant points}\label{alg_trisecant}
\begin{algorithmic}
\Procedure{Trisecant}{$\Rm$}\Comment{Find one trisecant in the Kummer variety of $X_{\Rm}$}
	\State $\Rm \gets \mathtt{siegeltrans}(\Rm)$ \Comment{Find Siegel's transform of $\Rm$}
	\State $\delta \gets 10^{-10}$ \Comment{Set precision $\delta$}
	\State $\ell \gets 0.1 $
	\State $\Delta \ell \gets 0.1$
	\State $\mathtt{max\_iter} \gets 100$ \Comment{Set maximum number of iterations}
	\State $ \Delta_N \gets [\quad ] $
	\While{$\ell\geq 0.5$}
	\State $X \gets \ell/2 \left(  e_{g-2} + \Rm e_{g-2}  \right) $
	\State $Y \gets \ell/2 \left( e_{g-1} + \Rm e_{g-1} \right) $
	\State $Z \gets \ell/2 \left( e_{g} + \Rm e_{g} \right)$
	\State $\mathbf{x} \gets [X(3:g);Y(2:g);Z(2:g)]$	\Comment{Initial vector}
	\State $V \gets [X(1:2); Y(1); Z(1)]$ \Comment{Fixed components $X_1$, $X_2$, $Y_1$, $Z_1$}
	\State $ [F,JF,\mathbf{x},\mathrm{res}, \Delta] \gets \mathtt{fay}(\mathbf{x},V,\Rm) $
	\State $ \pmb{\varepsilon} \gets - JF \backslash F$
	\State $\mathtt{iter} \gets 1$
	\While{$\Vert \pmb{\varepsilon} \Vert < \delta $} \Comment{\textit{Newton's iteration}}
	\State $\mathbf{x} \gets \mathbf{x} + \pmb{\varepsilon}$
	\State $ [F,JF,\mathbf{x},\mathrm{res}, \Delta] \gets \mathtt{fay}(\mathbf{x},W,\Rm) $
	\State $ \pmb{\varepsilon} \gets - JF \backslash F$
	\If{$ \mathrm{res} < \delta$ or $\mathtt{iter}>\mathtt{max\_iter}$}
	\State End Newton's iteration
	\EndIf
	\State $\mathtt{iter}\gets\mathtt{iter}+1$
	\EndWhile
	\State $\Delta_N  \gets [\Delta_N , \Delta] $ \Comment{Store the smallest $\Delta$ for each $\ell$}
	\State $\Delta_{\min} \gets \min(\Delta_N)$
	\If{$\Delta_{\min} < \delta$}
	\State End algorithm
	\EndIf
	\State $\ell \gets \ell + \Delta \ell$
    \EndWhile
	\State \textbf{return} $\mathrm{res}$, $\Delta_{\min} $
\EndProcedure
\end{algorithmic}
\end{algorithm}

The following is the routine computing the function $F:W\subset\C^{3g-4} \to\C^{2^g}$, where the fixed components $X_1,X_2,Y_1,Z_1$ are entered as parameters.
\begin{algorithm}[H]
\caption{Algorithm to compute the function $F$}
\begin{algorithmic}
\Procedure{fay}{$\mathbf{x},V,\Rm$}
	\State $X \gets [V(1:2); \mathbf{x}(1:g-2)]$
	\State $Y \gets [V(3); \mathbf{x}(g-1:2g-3)] $
	\State $Z \gets [V(4); \mathbf{x}(2g-2:2g-4)] $
	\State $X \gets [[X]] $	\Comment{Keep $X,Y,Z$ in the fundamental domain}
	\State $Y \gets [[Y]] $
	\State $Z \gets [[Z]] $	
	\State $\mathbf{x} \gets [X(3:g);Y(2:g);Z(2:g)]$
	\State $[\lambda_1,\nabla_Y \lambda_1, \nabla_Z \lambda_1] \gets \mathtt{lambdafun}(Y,Z,\Rm) $
	\State $[\lambda_2,\nabla_Z \lambda_2, \nabla_X \lambda_2] \gets \mathtt{lambdafun}(Z,X,\Rm) $
	\State $[\lambda_2,\nabla_X \lambda_3, \nabla_Y \lambda_3] \gets \mathtt{lambdafun}(X,Y,\Rm) $
	\State $[KX, JKX] \gets \mathtt{kummer}(X,\Rm)$
	\State $[KY, JKY] \gets \mathtt{kummer}(Y,\Rm)$
	\State $[KZ, JKZ] \gets \mathtt{kummer}(Z,\Rm)$
	\State $F \gets KX + \frac{\lambda_2}{\lambda_1} KY + \frac{\lambda_3}{\lambda_1}  KZ  $
	\State $J_X F \gets JKX + \frac{1}{\lambda_1} KY*\nabla_X \lambda_2 + \frac{1}{\lambda_1} KZ*\nabla_X \lambda_3 $ 
	\State $J_Y F \gets \frac{\lambda_2}{\lambda_1} JKY + \frac{1}{\lambda_1} KZ*\nabla_Y\lambda_3 - \frac{1}{(\lambda_1)^2} \left[ \lambda_2 KY + \lambda_3 KZ \right]*\nabla_Y \lambda_1 $
	\State $J_Z F \gets \frac{1}{\lambda_1} KY*\nabla_Z \lambda_2 + \frac{\lambda_3}{\lambda_1} JKZ - \frac{1}{(\lambda_1)^2} \left[ \lambda_2 KY + \lambda_3 KZ \right]*\nabla_Z\lambda_1 $
	\State $JF = [J_X F, J_Y F, J_Z F]$
	\State $JF(:,[1,2,g+1,2g+1]) \gets [\quad] $ \Comment{Deleting columns corresponding to $X_1,X_2,Y_1,Z_1$}
	\State $\mathrm{res} \gets \Vert F \Vert$ \Comment{Residual of $F$}
	\State $\Delta \gets \mathrm{min} \left( \mathrm{svd} (KX,KY,KZ) \right) $ \Comment{Linear dependence of $KX$ $KY$, $KZ$}
	\State \textbf{return} $F, JF,\mathbf{x}, \mathrm{res}, \Delta$
\EndProcedure
\end{algorithmic}
\end{algorithm}
Where $\mathtt{kummer}(X,\Rm)$ is the routine computing $\vec{\Theta} (X)$ (which is expressed as a column vector) and its Jacobian matrix $J\vec{\Theta} (X)$. The operation $*$ is the matrix multiplication.
\begin{algorithm}[H]
\caption{Algorithm to compute the Kummer map}
\begin{algorithmic}
\Procedure{kummer}{$X,\Rm$}
	\For{$I=1:2^g$}
	\State Set characteristic $\epsilon_I$	\Comment{$\epsilon_I \in \{0,1\}^g $ in column vector form}
	\State $factor \gets \exp \left( \frac{1}{2} \pi\I \langle \epsilon_I,\Rm \epsilon_I \rangle + 2\pi\I \langle \epsilon_I,X\rangle \right) $
	\State $[t, \nabla t] \gets \texttt{thetagrad}(2X+\Rm \epsilon_I,2\Rm,0) $ 
	\State $K(I) \gets factor \cdot t $
	\State $JK(I,:) \gets 2 \cdot factor\cdot (\nabla t + t \cdot \epsilon_I^T) $
	\EndFor
	\State \textbf{return} $K, JK$
\EndProcedure
\end{algorithmic}
\end{algorithm}
The routine $\mathtt{thetagrad}(X,\Rm,C)$ computes the multidimensional theta function $\Theta[C](X,\Rm)$ defined in \eqref{theta} with characteristic $C=\big[\begin{smallmatrix} \mathrm{p} \\\mathrm{q} \end{smallmatrix}\big]$, as well as its gradient $\nabla \Theta[C](X,\Rm)$ (which is expressed as a row vector) through \eqref{grad_theta}. We used the relation \eqref{thchar} to compute the even theta functions more efficiently in terms of the zero-characteristic theta function, since in practice we enter the array containing $\langle \mathrm{N},\mathbb{B} \mathrm{N}\rangle$ as an input, rather than computing $\langle ( \mathrm{N}+\epsilon_I/2),\mathbb{B}( \mathrm{N}+\epsilon_I/2)\rangle$ for every $\epsilon_I\in Z^g/2\Zg$, which would take up more memory space. Finally, the routine to compute the function $\lambda(a,b) := \Theta^*(a+b,\Rm) \Theta^*(a-b,\Rm)$ with and odd characteristic $2\delta_o\in \Z^{2g}/2\Z^{2g}$, as well as its gradients $\nabla_a \lambda(a,b)$, $\nabla_b \lambda(a,b)$ is the following:

\begin{algorithm}[H]
\caption{Algorithm to compute the lambda coefficients}
\begin{algorithmic}
\Procedure{lambdafun}{$a,b,\Rm$}
	\State $ \delta_o \gets [e_1/2; e_1/2] $ \Comment{Default odd characteristic}
	\State $[t_p, \nabla t_p] \gets \mathtt{thetagrad}(a-b,\Rm,\delta_o)$ 
	\State $[t_m, \nabla t_m] \gets \mathtt{thetagrad}(a+b,\Rm,\delta_o)$ 
	\State $\lambda \gets t_p \cdot t_m $
	\State $ \nabla_a \lambda \gets t_m \nabla t_p + t_p \nabla t_m $
	\State $ \nabla_b \lambda \gets t_m \nabla t_p - t_p \nabla t_m $
	\State \textbf{return} $\lambda,\nabla_a \lambda,\nabla_b \lambda$
\EndProcedure
\end{algorithmic}
\end{algorithm}
As before, since we use the same odd characteristic throughout the whole algorithm, we enter $\langle ( \mathrm{N}+e_1/2),\mathbb{B}( \mathrm{N}+e_1/2)\rangle$ as an input in the computation of the odd theta function $\Theta^*(\mathrm{z},\Rm)$ and its gradient.

\end{document}